\documentclass[10pt]{amsart}
\usepackage{amsfonts}
\usepackage{color}
\usepackage{cancel}

\theoremstyle{plain}
\newtheorem{theorem}{Theorem}[section]
\newtheorem{proposition}[theorem]{Proposition}
\newtheorem{lemma}[theorem]{Lemma}
\newtheorem{corollary}[theorem]{Corollary}

\theoremstyle{definition}
\newtheorem{definition}[theorem]{Definition}

\theoremstyle{remark}

\catcode`\@=11

\def\@ceqnnum{{\reset@font\rm (\tempequation)}}

\def\@ceqncr{{\ifnum0=`}\fi\@ifstar{\global\@eqpen\@M
    \@cyeqncr}{\global\@eqpen\interdisplaylinepenalty \@cyeqncr}}

\def\@cyeqncr{\@ifnextchar [{\@cxeqncr}{\@cxeqncr[\z@]}}

\def\@cxeqncr[#1]{\ifnum0=`{\fi}\@@ceqncr
   \noalign{\penalty\@eqpen\vskip\jot\vskip #1\relax}}

\def\@@ceqncr{\let\@tempa\relax
    \ifcase\@eqcnt \def\@tempa{& & &}\or \def\@tempa{& &}%
      \else \def\@tempa{&}\fi
     \@tempa \if@eqnsw\@ceqnnum\fi
     \global\@eqnswtrue\global\@eqcnt\z@\cr}

\def\clabel#1#2{\global\let\tempequation #2
   \@bsphack\if@filesw {\let\thepage\relax
   \def\protect{\noexpand\noexpand\noexpand}%
   \edef\@tempa{\write\@auxout{\string
      \newlabel{#1}{{#2}{\thepage}}}}%
   \expandafter}\@tempa
   \if@nobreak \ifvmode\nobreak\fi\fi\fi\@esphack}

\def\ceqnarray{
\global\@eqnswtrue\m@th
\global\@eqcnt\z@\tabskip\@centering\let\\\@ceqncr
$$\halign to\displaywidth\bgroup\@eqnsel\hskip\@centering
  $\displaystyle\tabskip\z@{{}##}$&\global\@eqcnt\@ne
  \hskip 2\arraycolsep \hfil${{}##}$\hfil
  &\global\@eqcnt\tw@ \hskip 2\arraycolsep
$\displaystyle\tabskip\z@{{}##}$\hfil
   \tabskip\@centering&\llap{##}\tabskip\z@\cr}

\def\endceqnarray{\@@ceqncr\egroup$$\global\@ignoretrue}

\catcode`\@=12

\long\def\salta#1{\relax}

\newcommand{\elle}[1]{L^{#1}(\Omega)} 
 
\newcommand{\loc}[1]{L^{#1}_{\rm loc}(\Omega)} 
\newcommand{\na}{\mathbb{N}}

\newcommand{\re}{\mathbb{R}}

\def\sob{W^{1,p}_{0}(\Omega)}

\def\gku{G_k (u)}

\def\lio{L^{\infty} (\Omega)}

\def\al{\alpha} 
\def\rn{\mathbb{R}^{N}} 
 
\def\D{\nabla} 
\def\sign{{\rm sign}\,} 
\def\vp{\varphi} 
\def\rife#1{(\ref{#1})} 
 
\def\eps{\varepsilon} 
\def\dive{\text{ \rm div }}

\def\dys{\displaystyle} 
 
\def\mis{\text{\rm{meas}}} 
\def\ene{L^p (0,T ; W^{1,p}_{0} ( \Omega ))} 
\def\eneloc{L^p (0,T ; W^{1,p}_{\rm loc} ( \rn ))}

\def\intbn{\displaystyle \int_{Q^T_n} }
\def\intbnp{\displaystyle \int_{Q^t_n} }

\def\intO{\dys \int_{\Omega}}

\newcommand{\locn}[1]{L^{#1}_{\rm loc} (\rn)}

\def\bc{\begin{cases}} 
\def\ec{\end{cases}} 
\def\be{\begin{equation}} 
\def\ee{\end{equation}} 
\def\vpla{\varphi_\lambda} 
\def\ga{\gamma}

\def\ga{\gamma}
\def\la{\lambda}

\def\t1pn{\mathcal T^{1,p}_{\rm loc} ((0,T)\times\rn)}

\def\ol{\overline}
\def\vare{\varepsilon}

\def\iakr{\int_{A_{k,R+\rho}}}
\def\iakrt{\int_{0}^{t_{1}}\int_{A_{k,R+\rho}}}

\title[Local estimates for parabolic problems with gradient terms]{Local estimates for  parabolic equations with nonlinear gradient terms}

\author[T. Leonori]{Tommaso Leonori}
\address{Tommaso Leonori, Departamento de An\'alisis Matem\'atico,  
Universidad de Granada, 18071   Granada, Spain}{} 
\email{leonori@ugr.es}
\thanks{The first author has been partially supported by MICINN Ministerio de Ciencia e Innovaci\'on (Spain) MTM2009-10878 and Junta de Andaluc\'{\i}a FQM-116.}

\author[F. Petitta]{Francesco Petitta}\address{Francesco Petitta, Departamento de An\'alisis Matem\'atico,
Universitat de Valencia,
C/ Dr. Moliner 50, 46100,
Burjassot, Valencia, Spain}\email{francesco.petitta@uv.es}
\thanks{ The second author has been partially supported by the PNPGC project, references MTM2008--03176.}

\begin{document}

\begin{abstract}
In this paper we deal with local estimates for parabolic problems in $\rn$ with absorbing first order terms, whose model is 
$$
\left\{
\begin{array}{ll}
u_t- \Delta u +u |\D u|^q  = f(t,x) \quad &\mbox{in}\,  (0,T) \times \rn\,, \\[1.5 ex]
u(0,x)= u_0 (x) &\mbox{in}\,   \rn \,,\quad 
\end{array}\right.
$$
where $T>0$, $N\geq 2$,  $1<q\leq 2$, $f(t,x)\in L^1(  0,T; L^1_{\rm loc} (\rn))$ and $u_0\in L^1_{\rm loc} (\rn)$.
\end{abstract}

\maketitle
\section{Introduction}

In this paper we deal with local estimates for parabolic problems in $\rn$, {$N\geq 2$},  with absorbing first order  lower order  terms. In particular, our main goal concerns the proof of the existence of a solution for  Cauchy problems whose model is 
\begin{equation}\label{model}
\left\{
\begin{array}{ll}
u_t- \Delta_p u +u |\D u|^q  = f(t,x) \quad &\mbox{in}\, (0,T) \times \rn\,, \\[1.5 ex]
u(0,x)= u_0 (x) &\mbox{in}\,   \rn \,,
\end{array}\right.
\end{equation}
where $T>0$, $p>1$, $p-1<q\leq p$, $f(t,x)\in L^1(  0,T ; L^1_{\rm loc} (\rn))$ and $u_0\in L^1_{\rm loc} (\rn)$, without any prescribed behavior of the solutions at infinity.  Here $\Delta_p u \equiv {\rm div}(|\nabla u|^{p-2}\nabla u )$ is the usual $p$-{laplace} operator. 
 
Such a problem is obviously strictly related to the possibility of proving estimates for the solutions that are independent  of   {their} behavior at infinity; this is a  peculiarity of nonlinear equations with  strong  \emph{absorption} lower order terms. 
{For instance, if    we consider the heat equation  in the whole space $\rn$, it is well known that the solution is explicit and turns out to be the convolution of the data with the heat kernel. 
Thus it is clear that, roughly speaking, $\forall (t,x)\in (0,T)\times \rn$, any  \emph{a priori} estimates on the solution  also depend on what happens  \emph{far away}  from $(t,x)$, which means that local estimates do not, in general,  hold true. 
}
\medskip

 Therefore,    the presence of the  the absorption lower order term,  and its   regularizing effect,    is crucial in order to prove local estimates; actually it plays the role of a barrier. 

\medskip

  If such a term does not depend on the gradient, i.e. for problems of the type 
\begin{equation}\label{semip}
\left\{
\begin{array}{ll}
u_t- \Delta_p u + b( u ) = f(t,x) \quad &\mbox{in}\, (0,T) \times \rn \\[1.5 ex]
u(0,x)= u_0 (x) &\mbox{in}\,   \rn \,,
\end{array}\right.
\end{equation}
 with $f(t,x)\in L^1(  0,T ; L^1_{\rm loc} (\rn))$, $b(\cdot) \in C^0 (\re)$,  and $u_0\in L^1_{\rm loc} (\rn)$, the existence (and regularity) of distributional solutions has been  investigated in \cite{bgv2} and more recently in \cite{lp}.

The main assumptions on the nonlinearity $b(u)$ are a sign condition  (namely  $ b( s )s \geq 0$, $\forall s\in \re$) and  an hypothesis on the behavior at  infinity: 
\be\label{leopel}
\int^{\infty} \frac{ds}{(b( s ) s)^{\frac1p}}<\infty\,.
\ee
We recall that if $p=2$ (and if $b$ is increasing at least at infinity) \rife{leopel} is equivalent to the well-known  Keller-Osserman condition. Such a condition has been introduced in the {pioneering} papers \cite{K} and \cite{o} in order to prove a local (uniform) bound for any subsolution of the  nonlinear  elliptic equation 
\be\label{lsell}
- \Delta u + b( u ) = f\quad \mbox{in}\,\, \Omega \,,
\ee
where $\Omega\subset\rn$ is   bounded and $f\in \lio$. 
This tool is strictly related to the possibility of constructing solutions that blow-up at the boundary (the so called {\it large solutions}); since the literature on this topic is huge,  we only mention, among the others, \cite{BaM} and \cite{Ve}. 
We want to stress that \rife{leopel} (as well as the Keller-Osserman condition) is the necessary and sufficient condition to have a global solution for  the ordinary differential equation associated to the (elliptic) equation, namely 
$$
\begin{cases}
(|v'|^{p-2} v')'= b(v)\quad \mbox{in }\ (0,\infty)\,,\\
v(0)=+\infty\,.
\end{cases}
$$
As we have already mentioned,   \emph{local estimates} are crucial in the study of large solutions for  \rife{lsell} and  such a problem turns out to be strictly related to the study of the same  problem in the whole $\rn$ without conditions at infinity (see \cite{Br}, \cite{BGV} and \cite{fab}).

\medskip

On the other hand,  nonlinear  equations with absorption gradient terms have been studied since many years. It has been recently  investigated (see \cite{tom}) the problem of both existence of solutions in the whole space ({\it entire solutions})  and existence of  large solutions associated to    nonlinear elliptic equations whose model is  
$$
-\Delta_p u + u+ u|\D u|^q = f(x) \quad \mbox{in}\, \Omega\,,
$$
where $\Omega\subseteq \rn$, $p>1$, $p-1<q\leq p$  and  $f(x)$ is   a singular datum (say $L^1 (\Omega)$). 

\medskip

The purpose of this paper is twofold: on one side we want to  extend the results of \cite{lp} to nonlinear problems with lower order terms that depend also on the gradient, namely of the type
\be\label{equa}
u_t-\Delta_p u +  u|\D u|^q = f(t,x) \quad \mbox{in}\, (0,T)\times \Omega\,,
\ee
for $\Omega =\rn$,  equipped with an initial datum $u_0 (x) \in \loc1$ and without any prescribed behavior at infinity for solutions.
 On the other hand, since we deal with local estimates, our aim is  to show that we can construct solutions that assume, in a suitable sense, the value \lq\lq $+\infty$\rq\rq{ }at $(0,T)\times \partial \Omega$, if $\Omega$ is a bounded domain in $\rn$.

\medskip
In order to prove such kind of results, we have to face several difficulties: first of all we have to {give a consistent definition of solution}.
In this framework several notions of solutions have been considered and, according with the definition introduced in  \cite{bm} for parabolic equations with $L^1$ data (see also  \cite{dpp}, \cite{dpr} and \cite{fra} for further generalizations to measure data), we use a renormalized formulation. 

In fact, such a notion  of  {  solution} turns out to be stronger than the distributional one, and it is very useful in order to face problems that involve singular data. 
Indeed, the peculiarity of such solutions is that they  do not have  local finite energy, but a priori estimates show that  their \emph{truncations}  belong to the energy space; the idea is then  to focus the attention to  a family of problems     solved by  such   truncations (see Definition \ref{def} below).
 
\medskip

The strategy  we use in order to prove the existence results relies on the combination of both  local estimates and  local compactness results in suitable Sobolev spaces for solutions of certain approximating problems.  It is  clear that, 
since  the formulation will involve cut--off functions, it  will be independent of the behavior at infinity, if we are interested in the problem in the whole space. 
We  also deal with the existence of large solutions for \rife{equa}  in the case of $ \Omega$ being a bounded domain. 
As far as large solutions are concerned,  we have the additional problem of defining the \emph{explosive} condition on $\partial \Omega \times (0,T)$: since, a comparison principle does not hold true, we can not approach it by using a sub and supersolutions method. 
Moreover,  the renormalized solutions are not,  in general, continuous, so we have to define \emph{how the value $+\infty$ is achieved on the boundary} in a convenient way.  For this purpose we extend the definition given in \cite{tom} to the parabolic framework; roughly speaking, we say that  a solution attains the value $+\infty$ on the boundary if  its truncation at level $k$ has constant trace $k$ on the boundary, for any $k>0$.  Because of this fact, in such a case, the truncations of the solutions have also to satisfy a global energy estimate in the whole $(0,T)\times \Omega$.
\medskip

\setcounter{equation}{0}
\section{Assumptions and statement of the main results }

Let $\Omega$ be an open subset of $\rn$, $N\geq 2$, possibly $\rn$ itself, and let $T>0$.    Throughout the paper we will use the following notation: for any ball $B_r$ of radius  $r>0$ and  $\forall \tau\in (0,T]$,  $Q^\tau_r =   (0,\tau)\times B_r$, while  $Q^\tau_{\Omega} =   (0,\tau)\times \Omega$.  

\medskip
We consider the following equation 
\be\label{eq}
u_t - \dive a (t, x ,u,\D u)   +g(t,x,u,\D u) = f(t,x) \quad \mbox{in}\ \,Q^T_{\Omega}\,, 
\ee
where $f(t,x)\in L^1(   0,T ; L^1_{\rm loc} (\Omega))$ while   $a (t, x ,s,\varsigma):[0,T] \times \Omega \times \re\times\rn \to \rn$ is a  Carath\'eodory function
  such that:
\begin{equation}\label{a1}
\exists \al> 0\,:\quad a (t, x ,s,\varsigma)\cdot \varsigma \geq \al |\varsigma|^p \,,
\end{equation}
\begin{equation}\label{a2}
\exists \beta > 0\,:\quad |a (t, x ,s,\varsigma) | \leq \beta |\varsigma|^{p-1} \,,
\end{equation}
and 
\begin{equation}\label{a3}
(a (t, x ,s,\varsigma)- a (t, x ,s,\eta))\cdot (\varsigma-\eta) >0\,,
\end{equation}
 for a.e. $(t,x) \in Q^T_\Omega$,    $\forall s \in \re $, and $ \forall \varsigma \, , \eta \in \rn\,$ such that $\  \eta \neq \varsigma$ and $p>1$. 
 
 \medskip 
 
Under the above assumptions    div $a (t,x,u,\D u)$  turns out to be a Leray-Lions type operator, acting from $L^p (0,T ; W_0^{1,p} (\Omega))$ into its dual  (see \cite{ll}).
\medskip

As far as the lower order term is concerned we suppose that   $g (t, x ,s,\varsigma):[0,T] \times \Omega \times \re\times\rn \to \re$ is a  Carath\'eodory function that satisfies the following;  $\exists L>0 $ such that 
\begin{equation}\label{g1}
g(t,x,s, \varsigma ) s \geq 0 ,\ \forall |s|\geq L\,, \, \forall \varsigma \in \rn \,, \mbox{ for a.e.  }  \, (t,x) \in Q^T_\Omega\,, 
\end{equation}
\begin{equation}\label{g3}
\begin{array}{c}
\quad \forall k>0,  \dys \sup_{|s|\leq k} |g(t,x,s,\varsigma)| \leq |g_k |+ \ga_k |\varsigma|^p, \, \,  \forall \varsigma \in \rn, \mbox{for a.e.} (t,x) \in Q^T_\Omega, \\
\ga_k>0, \quad    g_k (t,x) \in L^1(   0,T ; L^1_{\rm loc} (\rn)) \,, 
\end{array}
\end{equation}
and 
\begin{equation}\label{g2}
\dys   |g(t,x,s, \varsigma)| \geq  h(|\varsigma|^{p-1}) \,,\,\, \forall |s| \geq L,  \forall \varsigma \in \rn  \,, \mbox{ for a.e. }  (t,x) \in Q^T_\Omega\,, 
\end{equation} 
where $h$ is a
 {positive} $C^2 (\re^+)$ convex function $h$ such that $h(0)=0$, 
 and   the following conditions at infinity are satisfied:
\begin{equation}\label{g4}
\int^{\infty} \frac{d\tau}{h (\tau)} <\infty\, 
\quad \mbox{
and 
}\quad 
 \limsup_{\tau\to \infty} \frac{\tau^2 h'' (\tau)}{h'(\tau)\tau-h(\tau)}<\infty \,.
\ee

 Some comments {about} these assumptions are in order to be given.  Note that the absorption nature of the nonlinear lower order term depends on the sign condition \rife{g1}, while \rife{g3}   is known, in literature, as  \emph{natural growth condition}. 
We observe that condition \rife{g2} is a growth bound \emph{from below} for the lower order term with respect to $\varsigma$ at infinity. This assumption is crucial, as it can be noticed in the proof of   our main   results, and, in particular, it plays  a fundamental role in the  construction of suitable cut-off functions needed   to prove the local estimates we are interested in.

We remark  that the first condition  in \rife{g4} corresponds to the already mentioned   Keller-Osserman assumption for equation \rife{lsell}. In fact, in the same spirit of the stationary case,  it has to be  imposed in order to prove the existence of a solution for the ordinary differential equation associated to \rife{eq}. For instance,  a solution for the problem 
$$
\bc
\Big( |v'(s)|^{p-2 } v' (s)\Big)'= h \left(|v'(s)|^{p-1}\right) \qquad  \mbox{in}\, (0,+\infty)\,,\\
\dys \lim_{s\to0^+} v(s)=+\infty \,, 
\ec
$$
is well defined if and only if the first assumption in \rife{g4} holds true.   Finally, the second assumption in \rife{g4} is technical and we expect that it could  be removed.

\medskip 

\noindent Let us recall  the standard notation for truncations, i.e. $T_k (s)= \max \{ -k \min \{k, s\}\}$. 
We will, in general, handle with measurable functions whose truncations (locally) belong to the energy space $L^p (0,T;W^{1,p}_{\rm loc}(\rn))$. 
To do that it is useful to recall the notion of \emph{generalized gradient} whose main feature is contained in the following  result (see the proof in \cite{b6}, Lemma 2.1).

\begin{lemma}\label{bi6} Let $\Omega\subseteq \rn$, $N\geq 2$, and 
let  $w(t,x)$ be a measurable a.e. finite function such that   $T_k (w)\in L^1(0,T;W_{\rm loc}^{1,1}(\Omega))$, $\forall k>0$. 
Then there exists a measurable function $v:Q^T_\Omega\mapsto\rn$ such that $\forall k>0$ and for a.e. $(t,x)\in Q^T_\Omega$, 
\begin{equation}\label{benilan}
\nabla T_{k}(w)= v\chi_{\{|w|\leq k\}}\,.
\end{equation}
\end{lemma}
 Even if not explicitly stated, we will  made use of this notion throughout the paper.  
Anyway, we recall that, if $w\in L^1 (0,T;W_{\rm loc}^{1,1}(\Omega))$, then the generalized gradient coincides with the classical distributional one.
\medskip

\noindent Let us introduce the following definition of \emph{renormalized solution} which is the natural extension of the classical one (see {\cite{bm}}, \cite{dpp} and \cite{fra}).

\begin{definition}\label{def}
 We   say that   a measurable function $u (t,x) \in L^{\infty}( 0,T ; L^1_{\rm loc} (\Omega)) $ such that $\forall k>0$, $T_k (u ) \in L^p (0,T ; W^{1,p}_{\rm loc} (\Omega))$ is a \emph{renormalized solution} for equation \rife{eq}, if  $a (t, x ,u, \D u )\in (L^1( 0,T ; L^1_{\rm loc} (\Omega)))^N$, both  $f(t,x)$ and  $g(t,x,u,\D u) $ belong to $ L^1( 0,T ; L^1_{\rm loc} (\Omega))$, and the following identity holds true:
\be \label{renorm}
\begin{array}{c}
\dys 
-\int_{\Omega} S (u_0) \psi (x,0)
- \int_0^T  \langle S (u)  \,,\, \psi_t \rangle  \\[2.0 ex]
\dys + \int_{Q_\Omega^T} a (t, x ,u, \D u )\cdot  \D u S'' (u) \psi 
+ \int_{Q_\Omega^T}  a (t, x ,u, \D u )\cdot   \D \psi S' (u)  \\[2.0 ex]
\dys +\int_{Q_\Omega^T}  g(t,x,u ,\D u ) S '(u) \psi
=\int_{Q_\Omega^T}  f (t,x)S' (u) \psi\,,
\end{array}
\ee
$\forall \psi \in C^{1}_0 ( [0,T)\times \Omega)$, and for any 
$S(\tau)\in W^{2,\infty} (\re)$ such that  $S'(\tau)$ is compactly supported on $\re$.
Moreover,  
\be\label{fettine}
\lim_{l\to +\infty}  \int_{ Q_\Omega^T   \cap \{ l\leq |u|\leq l+1\}} a (t, x ,u,\D u)\cdot \D u  \,  \Psi =0 \,, \quad \forall \Psi \in  C^{0}_0 ( [0,T) \times \Omega )\,.
\ee
\end{definition}

Note that the regularity required for the solution is such that any term in \rife{renorm} makes sense. In fact the above definition is nothing  but  equation \rife{eq} formally multiplied by $S'(u)\psi $ and integrated on the cylinder $Q^T_{\Omega}$. The fact that $S'$ is compactly supported ensures that all but the first two terms in \rife{renorm} involve only a truncature of $u$. 
Condition \rife{fettine} is necessary to recover  a uniform  information on $u$ on the set where it is \emph{large}.

Finally, some comments regarding the initial datum are in order to be given: a priori, we are  not in the position to apply  Theorem 1.1 in   \cite{via} in order to deduce that $u \in C^0 ([0,T] ; L^1_{\rm loc} (\Omega))$, since we have not imposed any regularity on $u_t$. Anyway, this result can be applied to  $S(u)$, for any $S$ as above, since, by the equation,  the distributional time-derivative $S(u)_t$  turns out to belong to  $L^1(0,T; L^1 (\omega)) + L^{p'} (0,T;W^{-1,p'}(\omega))$, for any $\omega\subset\subset \Omega$.  
Actually, this is enough in order to give sense to the formulation, but does not imply any information about the continuity of $u$. 

\medskip

Here we state our existence result concerning entire solutions.

\begin{theorem}
\label{esistenza}
Assume that $a(t,x,s,\varsigma)$ and $g(t,x,s,\varsigma)$ satisfy  
\rife{a1}--\rife{a3} and $\,$ \rife{g4}--\rife{g2}, respectively. Then for any $f\in L^1 (0,T;L^1_{\rm loc} (\rn)) $ and for any $u_0\in\locn1$ there exists  a renormalized solution $u$ of the  Cauchy problem 
\begin{equation}\label{main}
\left\{
\begin{array}{ll}
u_t - \dive a (t, x ,u,\D u)   +g(t,x,u,\D u) = f(t,x) \quad &\mbox{in}\,(0,T) \times   \rn\\[1.5 ex]
u(0,x)= u_0 (x) &\mbox{in}\, \,  \rn \,.
\end{array}\right.
\end{equation}
Moreover $u \in C^0([0,T]; L^{1}_{loc}(\rn))$.
\end{theorem}

\medskip

As a consequence of the local estimates proved in the previous result, we are able to show the existence of a large solution for the boundary value problem associated to equation \rife{eq}. More precisely, let $\Omega \subset \rn$ be a bounded open subset of $\rn$, $N\geq 2$. We consider the following problem:

\begin{equation}\label{ls}
\left\{
\begin{array}{ll}
u_t - \dive a (t, x ,u,\D u)   +g(t,x,u,\D u) = f(t,x) \quad &\mbox{in}\,(0,T) \times   \Omega\\[1.5 ex]
u(t,x)= +\infty &\mbox{on}\ \,\partial_{\mathcal{P}} Q_{\Omega}^T \\[1.5 ex]
u(0,x)= u_0 (x) &\mbox{in}\,   \Omega \,.
\end{array}\right.
\end{equation}
where $\partial_{\mathcal{P}} Q_{\Omega}^T$ {denotes} the parabolic  vertical   boundary $(0,T)\times \partial{\Omega}$.

In the sequel we will need a suitable  version of \rife{g3} adapted for this context, namely  
\begin{equation}\label{g31}
\begin{array}{c}
 \forall k>0  
\quad \dys \sup_{|s|\leq k} |g(t,x,s,\varsigma)| \leq |g_k (t,x)|+ \ga_k |\varsigma|^p  \, ,\\
\ga_k>0, \quad   g_k (t,x) \in L^1(  Q_{\Omega}^T ) \,.
\end{array}
\end{equation}

Let us also specialize the definition of renormalized solution to this particular  boundary value problem. 
 
To our knowledge large solutions for parabolic equations have been investigated, for semilinear equations, in \cite{bd} and \cite{nv}. However,  for such class of equations, solutions are continuous, so that the \emph{explosive} condition makes sense pointwise.

For our purpose,   we need to reformulate  this condition in a suitable weak sense adapted to our framework. 
  More precisely, the   value \lq\lq$u=+\infty$\rq\rq at $\partial \Omega$ is assumed through a condition on  the trace of $T_k (u)$.

\begin{definition}\label{defls}
Let $\Omega$ be a bounded open subset of $\rn$, $N\geq 2$. For any $f(t,x)\in L^1(   0,T ; L^1_{\rm loc} (\Omega))$, we define a  {\it renormalized large solution} for problem \rife{ls} to be a measurable function $u (t,x) $ such that $T_k (u ) \in L^p (0,T ;  W^{1,p}(\Omega))$, $a (t, x ,u, \D u )\in (L^1( 0,T ; L^1_{\rm loc} (\Omega)))^N$, $g(t,x,u_n,\D u_n)\in L^1(   0,T ; L^1_{\rm loc} (\Omega))$  and it satisfies both \rife{renorm} and \rife{fettine}. Moreover the boundary condition is assumed in the following sense:
\be\label{infty}
k-T_k (u ) \in L^p (0,T ;  W^{1,p}_{0} (\Omega))\qquad \forall k>0\,.
\ee
\end{definition}

Our result concerning the  existence of a  large solution is the following one. 

\begin{theorem}
\label{exls}
Assume that $a(t,x,s,\varsigma)$ and $g(t,x,s,\varsigma)$ satisfy  \rife{a1}--\rife{a3} and   $\,$ \rife{g1}, \rife{g2}, \rife{g4} and \rife{g31}.
Then for any $f\in L^1(   0,T ; L^1_{\rm loc} (\Omega))$ such that $f^-\in L^1(  Q_{\Omega}^T)$ and for any $u_0\in L^1_{\rm loc} (\Omega)$ such that  $u_0^- \in L^1 (\Omega)$ there exists  a (renormalized) large solution $u \in C^0([0,T]; L^{1}_{loc}(\Omega))$ of   problem \rife{ls}.
\end{theorem}

We are also interested in some regularity properties  for the renormalized solutions of \rife{eq} both if $\Omega$ is bounded and if $\Omega =\rn$. Thus, let us introduce, for any $0<q<\infty$,  the Marcinkiewicz space $M^{q}(Q_\Omega^T)$ as the space of all  measurable functions $f$  such that there exists $c>0$, with 
 $$
 \text{meas}\{(t,x)\in Q_\Omega^T: |f(t,x)|\geq k\}\leq \frac{c}{k^q }, 
 $$
for every positive $k$ endowed with the seminorm
$$
\|f\|_{M^{q}(Q_\Omega^T)}=\inf \left\{ c>0: \text{meas}\{(t,x): |f(t,x)|\geq k\}\leq \left(\frac{c}{k }\right)^q \right\}\, .
$$
Let us recall that,  if $\Omega$  is bounded, then for $q > 1$ we have the following continuous embeddings
\be\label{embe'}
L^{q}(Q_\Omega^T)\hookrightarrow M^{q}(Q_\Omega^T) \hookrightarrow L^{q-\vare}(Q_\Omega^T),
\ee
for every $\vare\in (0,q-1]$.

We stress that from the definition of renormalized solution we can not, a priori, deduce neither any summability properties nor that $u$ has some continuity property in time for $p$ small.
 However, the following  result holds. 

{
\begin{proposition}\label{marc}
Any renormalized  solution of \rife{eq} with initial datum $u_0\in L^1_{\rm loc}(\Omega)$ satisfies the following estimates:
$$
\dys \|u\|_{M_{\rm loc}^{p-1 +\frac{p}{N}}(Q_\Omega^T)}\leq c_1 \ \ \ \text{and}\ \ \ \|\nabla u\|_{M_{\rm loc}^{ p-\frac{N}{N+1}}(Q_\Omega^T)}\leq c_2,
$$
where $c_1$ and $c_2$ are positive constants only depending on $u_0, f, N, R, T$ and  $p$. Moreover if $p>2-\frac{1}{N+1}$ then  $u \in C^0 ([0,T];L^1_{\rm loc} (\Omega))$. 
\end{proposition}
}
 
As already mentioned,  for $1<p\leq 2-\frac{1}{N+1}$ the continuity with values in $L^1_{\rm loc}$ can not be deduced \emph{a priori} by embedding theorems since the gradient does not in general  belong to any Lebesgue space. 
 However  the definition of renormalized solutions we gave above is not affected since, 
in order to give sense to the formulation, 
it is only required to $S (u)$   to admit a  trace at $t=0$.  Nevertheless, as stated in Theorem \ref{esistenza}, we shall see that the solution we have found satisfies this condition for any $p>1$. 
\medskip

We finally want to investigate how the local summability of the datum $f(t,x)$ influences the local summability of the renormalized solutions. In particular we will show that the regularity of the solutions is, locally, the same of the solutions of  equation \rife{eq} with  $g\equiv 0$  and equipped with homogeneous Dirichlet conditions at the vertical boundary. 

The techniques we use are nowadays classic and follow, for instance \cite{gm} and \cite{bpp}. However, since a localization is needed, the role of the lower order term (and in particular the growth condition \rife{g4}) is crucial. Actually we will only sketch the proof of such result, underlining the main differences with the cases  treated both  in \cite{gm} and in \cite{bpp}.  
{
\begin{theorem} \label{esse}
 Suppose $2-\frac{1}{N+1}<p<N$, $q>1$, $m>1$, and suppose that
 $f(t,x)$ belongs to $L^{m}(0,T; L^{q}_{\rm loc} (\Omega))$.
 Then for any  renormalized
solution  of \rife{eq}  there exists $C_0$ (depending on $u_0$, $f$, $N$, $\Omega$ and $T$) such that, 
\begin{enumerate}
 \item[$(i)$] 
if 
\begin{equation}\label{327}
1<\frac{1}{m} +\frac{N}{pq}\leq 1+\frac{N}{pm},
\end{equation}
for any initial datum $ u_0 \in L^{Nq\frac{ p-2+m'}{Nm'-pq}}_{\rm loc} (\Omega)$,
then 
 $$
\|u\|_{L^s (0,T; L^{s}_{\rm loc}(\Omega))}\leq C_0 \,,\quad  \mbox{where}\,\  s= \frac{mq(N+p)+N(p-2)(q(m-1)+m)}{mN-pq(m-1)};
$$
\noindent moreover 
$$\begin{array}{c}
\|u\|_{L^{s_1} (0,T; L^{s_2}_{\rm loc}(\Omega))}\leq C_0 \,,\quad  \mbox{where}\,\  s_1=m' s_0\,,\, s_2 = q's_0 
\\[2.0 ex] \mbox{and} \,\  \dys s_0 =\frac{mq(q-1)+q(m-1)[p(N+1)-2N]}{mN-pq(m-1)}.
 \end{array}
$$ 
 \item[$(ii)$]  
If \begin{equation}\label{327n}
\frac{1}{m} +\frac{N}{pq}>1+\frac{N}{pm}.
\end{equation}
for any initial datum $ u_0 \in L^{\frac{N(q-1)(p-1)+N-pq}{N-pq}}_{\rm loc} (\Omega)$,
then 
 $$
\|u\|_{L^s (0,T; L^{s}_{\rm loc}(\Omega))}\leq C_0 \,,\quad  \mbox{where}\,\  s= \frac{[N(p-1) (q-1) +N -pq](N+p)}{N(N-pq)}+p-2;
$$
moreover 
$$
\|u\|_{L^{s_1} (0,T; L^{s_2}_{\rm loc}(\Omega))}\leq C_0 \,,\,   \mbox{where}\, \ s_1=m' s_0\,,\, s_2 = q's_0 \quad \mbox{and} \, 
 s_0 =\frac{N(q-1)(p-1)}{N-pq}.
$$ 
 \item[$(iii)$]  
If 
\begin{equation}\label{326}
\frac{1}{m} +\frac{N}{pq}<1\,,
\end{equation}
for any initial datum $ u_0 \in L^{\infty}_{\rm loc} (\Omega)$,
then 
 $
\|u\|_{L^{\infty} (0,T; L^{\infty}_{\rm loc}(\Omega))}\leq C_0.
$
 \end{enumerate}
\end{theorem}
}

Let us only  notice that, as a typical smoothing effect for parabolic problems, the summability of the solution in space (e.g.  $s_2$ in Theorem \ref{esse}) is always greater than the summability of the initial datum.

\subsection*{Notation} Define $\vp_{\la} (s) = s e^{ \la s^2}$; we recall  that  $\vp_{\la} (s)$ enjoys the  following useful property: 
\begin{equation} \label{vpla} 
\forall a > 0\,,\,b >0\,,\, \forall  \la > \frac{b^2}{8 a^2}  \qquad a \vp_{\la} ' (s)  - b |\vp_{\la} (s)|  \geq 1\, ,\quad \forall s \in \re\,.
\end{equation}

We will also make use of the following functions related with the truncations:
\be\label{sgei}
S_j (\tau)=\int_0^\tau [1-T_1 (G_j (s))] ds\,,
\ee
and $G_k (s)= s-T_k (s)$.

\medskip

By $\langle\cdot , \cdot \rangle$  we mean the duality between suitable spaces in which function are involved. In particular we will consider both the duality between $W^{1,p}_0 (\Omega)$ and $W^{-1,p'} (\Omega)$ and  the duality between $W^{1,p}_0 (\Omega)\cap L^{\infty} (\Omega)$ and  $W^{-1,p'} (\Omega)+L^{1} (\Omega)$.

\medskip

Finally, we  use the following notation for sequences:
 $\eps (\sigma, n, \nu)$ will  indicate any quantity that vanishes as the parameters go to
their (obvious, if not explicitly stressed) limit point,  with the same order in which they appear,  that is, 
$$
\dys\lim_{\nu\rightarrow \infty} \limsup_{n\rightarrow +\infty}
\limsup_{\sigma \to \infty} |\eps(\sigma,n,\nu)|=0.
$$
We will also sometimes omit the dependence of $\vare(\cdot)$ on one or more of its arguments, when they are not present.
\medskip

\setcounter{equation}{0}
\section{Technical results}

In this section we collect some technical results that will be useful in the rest of the paper. The first one concerns  the construction of a  suitable family of functions.

\begin{proposition}\label{prop}
Let $h  : \re^+ \to \re^+$ be a $C^2$, convex  function, such that  $h (0) = 0$, and such that  \rife{g4}   holds. Then, for any $\delta>0$, there exists a constant   $C_0=C_0 ({\delta}) >0$ and a function $\sigma: [0,1]\to[0,1]$, $\sigma \in C^0 ([0,1])\cap C^1 ((0,1))$ with $\sigma (0) = \sigma '(0) = 0$, $\sigma (1) = 1$,  such that 
\be\label{lemmaa}
\quad \forall v>0,\qquad v \sigma' (s) \leq \delta h (v) \sigma (s) + C_{\delta}\,,\quad \forall s\in [0,1]\,.
\ee
\end{proposition}

Before giving the proof of Proposition \ref{prop}, we need to introduce another fundamental tool in our arguments, that is a  generalized Young inequality with the function $h$ which appears in \rife{g2}-\rife{g4}. In order to do that,   we have to introduce the Legendre transform  for $h$ together  with its properties which we will use in the sequel. 

\medskip

We recall that  $h$ is a $C^2$ increasing and convex function such that $h(0)=0$. Moreover by the convexity and since \rife{g4} holds (i.e.,   roughly speaking, $h$ is \emph{a bit more} than  superlinear at infinity) it follows that 
 $$\lim_{s\to \infty} h' (s)=+\infty\,.$$
Let us consider the  {\it Legendre transform } of $h$ defined by 
$$
h ^{\ast} (q)=\sup_{r \in \re} \left[q r-h (r)\right]
\,.$$
Here we recall the so called generalized Young inequality; namely, for any positive $z, w $, we have
\be\label{young}
wz\leq h (z)+ h^* (w)\,.
\ee
It is clear that $h^{*}$ is continuous, increasing and, since \rife{young} holds, superlinear at infinity. 
Consequently  ${h^{*}}^{-1}$ is well defined and moreover 
$$\lim_{q\to \infty}{h^{*}}^{-1}(q)=+\infty\,.$$ 
Moreover, since $h$ is smooth,  $\forall q>0$,  we have 
$$ h ^{\ast} (q)= q [(h')^{-1} (q)]-h ((h')^{-1} (q))\,,$$ so that  $$h ^{\ast} (h' (y) )= y h'(y)- h(y)\,,\ \ \forall\ y>0\,.$$

{The proof of } Proposition \ref{prop} is based on the possibility of constructing a solution of  a  suitable Cauchy problem, as stated in the following Lemma. 

\begin{lemma}\label{cutoff}
Let $h  : \re^+ \to \re^+$ be a $C^2$, convex  function, such that  $h (0) = 0$, and such that  \rife{g4}    holds. Then, for any $\delta > 0$ there exists $C_0 = C_0 (\delta) $ and a function $\sigma=\sigma_{\delta}: [0,1]\to[0,1]$, $\sigma \in C^0 ([0,1])\cap C^1( (0,1))$ 
solution of the problem 
\be\label{lemma}
\bc
\dys {\sigma}' (s)= \delta  {\sigma} (s) {h^*}^{-1}\left(\frac{C_0}{\delta  {\sigma} (s)} \right)\quad \mbox{in} \, \,(0,1)\,, \\[2.0 ex]
 {\sigma} (0)=0\,,\quad \sigma (s)>0\,.
\ec
\ee
Moreover
\be\label{lemma2}
\sigma (1) = 1\,,\quad \mbox{and}\quad \lim_{s\to 0^{+}} \sigma ' (s)=0\,.
\ee
\end{lemma}

\proof

Let us consider the family of functions $\sigma (s)$ defined by the implicit formula
\be\label{wd}
\dys 
\int_0^{\sigma(s)}  \frac{dt }{(h^*)^{-1} \left( \frac{\tau }{t}\right) t} = s\delta \,, \quad \forall \tau >0\,.
\ee
Our aim is to prove that $\sigma$ is well defined and that there exists a   value  $\tau=C_0$ such that \rife{lemma} and \rife{lemma2} hold true.

\medskip 
\noindent{\bf Step 1: Near 0.}
We want to prove that $\forall \tau >0$, $\sigma (s)$ is well defined in a neighborhood of $s=0$. Indeed, through the change of variable defined by the relationship 
{$h' (z)=  (h^*)^{-1} \left(\frac{\tau}{\delta \sigma(s)}\right)$}, and by the properties of $h$ and $h^*$ we have stated before,  it follows that
{
$$\dys 
\frac{1}{\delta}\int_0 \frac{dt }{(h^*)^{-1} \left( \frac{\tau }{t}\right) t
}  <+\infty
\quad 
\Leftrightarrow
\dys \quad \frac{1}{\delta}
 \int^{+\infty} \frac{z h'' (z)}{h'(z)[h'(z)z-h(z)]} dz <+\infty\,.
$$}
Recalling  \rife{g4} and  since $h'(z)z-h(z) > 0$, for any $z>0$,    there exists a constant $c_1$ such that 
$$
\dys 
 \int^{+\infty} \frac{1}{z h'(z)} \frac{z^2 h'' (z)}{h'(z)z-h(z)} dz \,\leq c_1 \dys 
 \int^{+\infty} \frac{dz }{h'(z) z} dz \leq  c_1 \dys 
 \int^{+\infty} \frac{dz }{h(z) } \,,
$$
where the last inequality holds since $h$ is convex and $h(0)=0$:  by \rife{g4} the last integral is finite and so $\sigma$ is well defined near $0$.

\medskip 
\noindent{\bf Step 2: The choice of $C_0$.}
It follows by Step 1, through the change   $\rho = \frac{\tau}{\delta t}$, $\forall \delta >0$,  that
$$
\lim_{\tau\to+\infty} \frac1{\delta} \int_{\frac{\tau}{\delta}}^{\infty} \frac{d\rho }{\rho (h^{*})^{-1} (\rho)} =0\,;
$$
on the other hand, since $(h^{*})^{-1} (0)=0$, 
$$
\lim_{\tau\to0} \frac1{\delta} \int_{\frac{\tau}{\delta}}^{\infty} \frac{d\rho }{\rho (h^{*})^{-1} (\rho)} =+\infty\,.
$$
Thus there exists $C_0 $ such that 
$$
\frac1\delta \int_{\frac{C_0}{\delta}}^{\infty} \frac{d\rho }{\rho (h^{*})^{-1} (\rho)} =\dys 
\int_0^{1}  \frac{dt }{(h^*)^{-1} \left( \frac{C_0 }{t}\right) t
} =1\,,
$$
which implies $\sigma(1)=1$.

\medskip 
\noindent{\bf Step 3: The limit of $\sigma'$.}
Recalling the definition of $\sigma'$ from \rife{lemma}, we want to prove 
\be\label{sarcazzo}
\lim_{s\to 0^+} \delta  {\sigma} (s) {h^*}^{-1}\left(\frac{C_0}{\delta  {\sigma} (s)} \right)=0\,.
\ee
This is equivalent to prove that 
$$
\lim_{\tau \to +\infty} \frac{h' (\tau)}{h^*(h' (\tau))}  = \lim_{\tau \to +\infty} \frac{h' (\tau)}{h' (\tau)\tau -h (\tau)} =0
$$
since $\tau$ is such that $(h^{\ast})^{-1}(\frac{C_0}{\delta \sigma (s)})= h' (\tau)$. Using that $h^*(h' (\tau))\to +\infty $ as $\tau$ diverges and by De l'Hopital rule we deduce that \rife{sarcazzo} holds, and so the Lemma is proved.\qed

\medskip

\proof[of Proposition \ref{prop}].  
Let $\sigma (s)$ be the function defined in Lemma \ref{cutoff}. Thus it is clear that inequality  \rife{lemmaa} is satisfied at $s=0$, and we can multiply and divide the left hand side by $\sigma (s)$;  using  \rife{young} we get
$$
\delta \sigma (s)  v \frac{ \sigma' (s)}{\delta \sigma (s)} \leq  
 \delta \sigma (s) h(v) +  \delta \sigma (s) h^* \left(\frac{ \sigma' (s)}{\delta \sigma (s)}\right)\,.
$$
Recalling that $\sigma$ is the solution of the Cauchy problem defined in \rife{lemma},  \rife{lemmaa} holds true.
\qed
\medskip

In the sequel we will  also handle with dualities involving the time derivatives of suitable functions; to this aim we will  
use the following Landes-type (see \cite{lan}) regularization result.

\begin{lemma}\label{lg}  Let $\Omega$ be an open bounded subset of  $\rn$, and  let $w\in L^p (0,T; \sob)$ and $w_0\in\elle{1}$. Then, for any $\nu>0$, there exists a function $\eta_\nu= \eta_\nu (w, w_0)\in L^p (0,T; \sob)$, such that 
$$
\frac{d}{dt}\eta_\nu  =\nu(w- \eta_\nu  ), 
$$
and  $ \eta_\nu (w, w_0 )(0,x)=\eta_{0,\nu}\in\elle{2}$, with 
$$
\eta_{0,\nu} \stackrel{\nu \to \infty}{\longrightarrow} w_0   \ \ \text{in}\ \ \ L^1(\Omega). 
$$
If furthermore  $w\in L^{\infty} (Q_{\Omega}^T)$, then 
\be\label{lim}
\dys \|\eta_\nu\|_{ L^{\infty} (Q)} \leq \|w\|_{ L^{\infty} (Q_{\Omega}^T)}\,.
\ee
Moreover, if $w_t = w^{(1)}+w^{(2)}\in L^1(Q_{\Omega}^T)+L^{p'}(0,T;W^{-1,p'}(\Omega))$,  then     $\frac{d}{dt}\eta_\nu$ { admits a decomposition of the form} $\frac{d}{dt}\eta_\nu=\rho_{\nu}^{(1)}+\rho_{\nu}^{(2)}$, with both 
$$
\rho_{\nu}^{(1)}\stackrel{\nu \to \infty}{\longrightarrow} w^{(1)}    \ \ \text{in}\ \ \ L^1(Q_{\Omega}^T)
$$
and
$$
\rho_{\nu}^{(2)}\stackrel{\nu \to \infty}{\longrightarrow} w^{(2)}    \ \ \text{in}\ \ \ L^{p'}(0,T;W^{-1,p'}(\Omega)).
$$
\end{lemma}

\proof See  \cite{gre}, Lemma 2.1.
\qed 
\medskip 

Here we state a useful result which allows us to handle functions that do not have  time derivatives belonging to the dual of the energy space $L^p(0,T;\sob)$. In fact it consists in a generalized integration by
parts formula,  whose proof can be found in \cite{dpr} (see also
\cite{cwi}).

\begin{lemma}\label{cw}
Let $\Omega$ be any domain in $\rn$, $N\geq 2$, and 
let $\phi :\re\to\re$ be a continuous piecewise $C^1$ function such that $\phi (0)=0$ and $\phi'$ has compact support; let us define $\Phi(s)=\int_{0}^{s} \phi (r)dr$. If $v\in L^p(0,T,W^{1,p}_0 (\Omega))$ is such that $v_t \in L^{p'}(0,T;W^{-1,p'}(\Omega))+ L^1(Q^T_\Omega)$ and if $\psi\in C^{\infty}(\overline{Q^T_\Omega})$, then we have
\begin{equation}\label{ibp2}
\int_{0}^{T}\langle v_t,\phi (v)\psi\rangle\ =\int_\Omega \Phi(v(T))\psi(T)\  -\int_\Omega \Phi(v(0))\psi(0)\ 
-\int_{Q^T_\Omega}\psi_t \ \Phi(v)\ .
\end{equation}
\end{lemma}

We observe that   $v_t \in L^{p'}(0,T;W^{-1,p'}(\Omega))+ L^1(Q^T_\Omega)$ implies
that there exist $\eta_1\in L^{p'}(0,T;W^{-1,p'}(\Omega))$ and $\eta_2
\in L^1 (Q^T_\Omega)$ such that $u_t =\eta_1 +\eta_2$. Even if   $\eta_1$ and
$\eta_2$ are not uniquely determined,   the integration by parts
formula turns out to be independent of the representation of $v_t$; moreover, 
according with the notation introduced before, $\langle\cdot,\cdot\rangle$ will indicate the duality  between
$L^{p'}(0,T;W^{-1,p'}(\Omega))+ L^1(Q_\Omega^T)$ and $L^{p}(0,T;W^{1,p}_0 (\Omega))\cap L^{\infty} (Q_\Omega^T)$.

\medskip

We also recall the following classical result due to Gagliardo and Nirenberg.

\begin{theorem}[Gagliardo-Nirenberg]\label{gnet}
Let $\Omega\subset \rn$, open and bounded, 
and let $v$ be a function in $W^{1,\mu}(\Omega)\cap\elle\lambda$ with $\mu \geq 1$, $\lambda \geq 1$. Then there exists a positive constant $C$, depending on $N$, $q$ and $\lambda$, such that
 $$
 \|v\|_{\elle\eta}\leq C\|\nabla v\|_{(\elle{\mu})^N}^{\theta}\|v\|^{1-\theta}_{\elle\lambda}\,,
 $$ 
 for every $\theta$ and $\eta$ satisfying
 $$
 0\leq \theta\leq 1,\ \ \ 1\leq\eta\leq +\infty,\ \  \ \dys\frac{1}{\eta}=\theta\left(\frac{1}{\mu}-\frac{1}{N}\right)+\frac{1-\theta}{\lambda}\,.
 $$
 \end{theorem}
 \proof
 See \cite{n}, Lecture II.
 \qed

 The following embedding results are consequences of the previous theorem. We will use them in the last section but  we give here their statement  for completeness.

 \begin{corollary}\label{cgnet}Let $v\in L^q(0,T; W^{1,q}_0 (\Omega))\cap L^{\infty}(0,T;\elle\gamma)$, with $q\geq 1$, $\gamma\geq 1$. Then $v\in L^\sigma(Q^T_{\Omega})$ with $\sigma=q\frac{N+\gamma}{N}$ and 
\begin{equation}\label{gnc}
\int_{Q^T_{\Omega}} |v|^\sigma \ dxdt\leq C \|v\|_{L^{\infty}(0,T;\elle\gamma)}^{\frac{\gamma q}{N}}\int_{Q^T_{\Omega}} |\nabla v |^q \ dxdt\,.
\end{equation}\end{corollary} 

\begin{corollary}\label{lemma310}
Let $\Omega\subset \rn$, open and bounded,  $\tau > 0, 1 < p < N$ and let further 
$w 
\in L^{\infty} (0, \tau ; \elle{p}) 
\cap L^p 
(0, \tau ; W^{1,p}_0 (\Omega)) $.
 Then there exists a positive constant K depending 
only on N and p such that
$$\begin{array}{l}
\dys\left[\int_{0}^{\tau}\left(\intO |w|^{\sigma}\right)^{\frac{\mu}{\sigma}}\right]^{\frac{p}{\mu}} 
 \dys  \leq  K \left(\sup_{t\in[0,\tau]}\intO |w|^p + \int_0^{\tau} \intO |\D w|^p
\right) 
\end{array}
$$ 
for all $\mu$ and $\sigma$ satisfying 
\begin{equation}\label{330}
p\leq \sigma\leq p^{\ast},\ \ \ p\leq\mu\leq\infty,\ \ \ \frac{N}{p\sigma}+\frac{1}{\mu}=\frac{N}{p^{2}}.
\end{equation}
\end{corollary}

We  also recall the interpolation inequality that we will use in the proof of Theorem \ref{esse}.
Assume that $z\in L^\infty(0,T;L^q (\Omega))\cap L^p(0,T; L^r (\Omega))$, $p,q,r\geq 1$. Thus   $z\in L^{\eta} (Q_\Omega^T)$ and  
\be\label{interpolo}
\begin{cases}
\begin{array}{c}
\|z\|_{L^{\eta}(Q)}\leq C \|z\|_{L^\infty(0,T;L^q (\Omega))}^{1-\theta}\|z\|_{L^p(0,T; L^r (\Omega))}^{\theta}\\[2.0 ex]
 \mbox{
with }\quad \frac{1}{\eta}=\frac{\theta}{r}+\frac{1-\theta}{q}, \ p\geq \theta \eta. 
\end{array}
\end{cases}
\ee

A useful application of Corollary \ref{cgnet} is the following.
\begin{proposition} \label{p-1}
Let $\Omega \subset \rn$ be a bounded domain and  $p>1$. 
Let $w\in L^\infty(0,T; L^1(\Omega))$ such that $T_k(w)\in L^p(0,T;W_{0}^{1,p}(\Omega))$, for any $k>0$. If $|\nabla w|^{p-1}\in L^1(Q_{\Omega}^T)$, then $w^{p-1}\in L^1(Q_{\Omega}^T)$.
\end{proposition}

\proof
We deal only with the case $p>2$, since for $p\leq2$ it is trivial. Since both $w\in L^\infty(0,T; L^1(\Omega))$ and $|\nabla w|^{p-1}\in L^1(Q_{\Omega}^T)$, we have that  $w\in L^1 (0,T;W^{1,1}_0 (\Omega))$. Then, we can apply Corollary \ref{cgnet} with $q=\gamma=1$ to obtain that $w\in L^{\frac{N+1}{N}}(Q^T_{\Omega})$. Now, if $p\leq1 + \frac{N+1}{N}$ we are finished, otherwise $w\in L^1 (0,T;W^{1,\frac{N+1}{N}}_0 (\Omega))$ and we apply again Corollary \ref{cgnet} with $\gamma=1$ and $q=\frac{N+1}{N}$ to conclude that $w\in L^{\left(\frac{N+1}{N}\right)^2 }(Q_{\Omega}^T)$. It is clear that, iterating this procedure, we get the result in a finite number of steps.
\qed

\medskip
 
The estimates contained in the following lemma are standard and turn out to coincide, for instance, with the one proved  in \cite{bg}  (see also Lemma 3.7 in  \cite{lp}).
{
\begin{lemma}\label{gajardo}
Let $\Omega \subset \rn$ be a bounded domain  and let $w\in L^p (0,T;W^{1,p} (\Omega)) \cap L^{\infty} (0,T;L^{1}  (\Omega))$,  $1<p<N$. Suppose moreover that there exists $C_0>0$ depending only  on $ N, \Omega, T$ and  $p$   such that 
$$
\int_{Q_{\Omega}^T} |\nabla T_k(w)|^p\leq C_0(k+1) \quad \mbox{and} \quad \int_0^T |w| \leq C_0\,,
\quad \forall k>0\,.
$$
 Then:
$$
\dys \|w\|_{M^{s_1}  (Q_{\Omega}^T)} \leq c_1 \ \ \ \text{and}\ \ \ \|\nabla w\|_{M^{ s_2}(Q_{\Omega}^T)} 
\leq c_2$$
where  
$$
\bc
s_1= \max\left\{1,p-1 +\frac{p}{N}\right\} \\[1.5 ex]
s_2= \max\left\{ \frac{p}{2}, p-\frac{N}{N+1}\right\}\,, 
\ec
$$
and $c_1$ and $c_2$ are positive constants only depending on $C_0, N, \Omega, T$ and  $p$.
\end{lemma}}

Finally let us state the following  classical result due to Stampacchia.

\begin{lemma} \label{stampa}\sl
Let $\zeta (j,\rho)\ :[0,+\infty)\times[0,R) $ be a function such
that $\zeta (\cdot,\rho) $ is nonincreasing and $\zeta (j,\cdot)
$ nondecreasing. Moreover, suppose that $\exists K_0 >0$, $\mu >
1,$ and $C,\nu,\ga >0$ such that
$$
\zeta (j,\rho) \leq C \frac{ \zeta(k,R)^\mu}{(j-k)^\nu
(R-\rho)^\ga}\quad \forall j>k>K_0 ,\, \forall  \rho \in
(0,R].
$$
Then for every $  \delta \in (0,1)$, there exists $ d>0$ such that:
$$
 \zeta (K_0+ d,(1-\delta)R ) =0,
$$
where
$$
\dys d^\nu = C'  2^{\frac{\mu(\nu+\ga)}{\mu-1}}  \frac{
\zeta(K_0,1)^{\mu-1}}{(1-\delta) },\qquad C'>0.
$$
\end{lemma}

\proof See \cite{s}. \qed

\setcounter{equation}{0}
\section{proofs of Theorem \ref{esistenza} and Theorem \ref{exls}}

Here, and throughout the paper,  we denote by $\xi =\xi^{\rho}_R (|x|)$ a $C^{1}_0 (\rn)$ function such that $\forall \rho>0$ 
\be\label{xi}
\left\{
\begin{array}{ll}
\xi \equiv 1 \quad &\mbox{ if   }\, |x|\leq R\\[1.5 ex]
0< \xi<1 \quad &\mbox{ if   }\, R< |x|< R+\rho \\[1.5 ex]
\xi \equiv 0 \quad & \mbox{ if   } \, |x|\geq R+\rho\,.
\end{array}
\right.
\ee
Note that this class of functions keeps its properties if composed with $C^1[0,1]$ functions $w:[0,1]\to [0,1]$ such that:
$$
w(0)=w'(0)=0, \quad w(1)=1, \quad  \mbox{and} \, \quad w(s)>0 \, \,\mbox{in} \, \,(0,1)\,.
$$
We will often choose a cut-off function $\xi=\sigma (\eta)$ where $\eta$ (and consequently $\xi$) satisfies \rife{xi} and $\sigma$ is the function that appears  in  Proposition  \ref{prop}, and thus such that \rife{lemmaa} holds true for a suitable choice of $\delta$.

\medskip
Now we can prove the first result concerning existence of solutions.
\medskip
\proof[of Theorem \ref{esistenza}]. 
Let   $B_n$ be the ball of radius $n$, centered at the origin (this  choice can be done without loss of generality) and let   $u_n$ be a weak solution of the following problem
\begin{equation}\label{app}
\left\{
\begin{array}{ll}
(u_n)_{t}- \dive a (t, x ,u_n,\D u_n)   +g(t,x,u_n,\D u_n) = f_n(t,x) \, &\mbox{in}\,Q^T_n  ,\\[1.5 ex]
 u_n (t,x)= 0 & \mbox{on}\,    \partial_{\mathcal{P}} Q_n^T , \\[1.5 ex]
u_n (0,x)= u_n^0 (x) &\mbox{in}\,   B_n ,
\end{array}\right.
\end{equation}
where $f_n (t,x)= T_n (f(t,x))$ and $ u_n^0 (x) =T_n (u_0 (x))$. 
Note that,  thanks to the result of \cite{do} (see also \cite{via}), there  exists (at least) a weak solution for \rife{app}, i.e. a function $u_n \in L^p (0,T ; W^{1, p}_0 (B_n))$ such that $(u_n)_t$  {belongs to} $L^{p'} (0,T ; W^{-1, p' } (B_n))$, $g(t,x,u_n,\D u_n)$  {belongs to} $L^1 ((0,T)\times B_n)$, and the following identity holds true 
\be\label{weakn}
\begin{array}{c}
\dys \int_0^T  \langle  (u_n)_t   \,,\, \psi \rangle
+ \intbn a (t, x ,u_n, \D u_n )\cdot  \D  \psi  \\
+\intbn g(t,x,u_n,\D u_n)    \psi
=\intbn f_n    \psi\,, 
\end{array}
\ee
$\forall \psi \in L^p (0,T ; W^{1, p}_0 (B_n))\cap L^{\infty} (Q_n^T)$.
\medskip

We will prove Theorem \ref{esistenza} by showing that the terms in \rife{app} are compact in suitable spaces. In order to do it, here and throughout the whole proof, we fix a ball $B_R$, centered at the origin, and we will prove suitable estimates for $u_n$ in $(0,T) \times B_R$.  Moreover, a weak solution on $Q_n^T$
 turns out to be obviously a weak solution in $Q_n^t$ for any $0<t<T$. Hence, with an abuse of notation, we will often refer to \rife{weakn} by tacitly understanding its counterpart on $Q_n^t$.   
 
\medskip  
\noindent {\bf Local estimates on truncations.}
For any $n\geq R+\rho$ (for any fixed $\rho>0$),  let us choose in \rife{weakn} $ \psi= \vpla{(T_k ( u_n))} \xi$,  where $\vpla (s)=se^{\la s^2}$ ($\lambda >0 $ will be fixed later), $k>L$,  and $\xi (x)$ is  a cut-off function such that \rife{xi} holds true  (we will often omit the dependence  on $x$).
Thus we have
\be \label{est1}
\begin{array}{c}
\dys \int_0^t
\langle (u_n)_t , \vpla ( T_k (u_n)) \xi \rangle
\\[2.0 ex]
+ \intbnp  a (t, x ,u_n, \D u_n )\cdot  \D  T_k (u_n)   \vpla' ( T_k (u_n))   \xi
\\[2.0 ex]
+ \intbnp  a (t, x ,u_n, \D u_n )\cdot  \D  \xi   \vpla  (T_k (u_n)) 
\\[2.0 ex]
+\intbnp  g(t,x,u_n,\D u_n)     \vpla ( T_k (u_n)) \xi
=\intbnp  f_n(t,x)     \vpla ( T_k (u_n)) \xi\,.
\end{array}
\ee
Since $\xi$ does not depend on time, using Lemma \ref{cw}, 
$$
\int_0^t
\langle  (u_n)_t , \vpla (T_k (u_n)) \xi \rangle= \int_{B_n } \Phi_{\lambda, k} ( u_n (t,x) ) \xi -  \int_{B_n } \Phi_{\lambda, k} (   u_n (x, 0)) \xi \,,
$$ where 

$$
\Phi_{\lambda,k} (s)= \int_0^s \vpla (T_k (\tau)) d\tau= 
\left\{
\begin{array}{ll}
\frac{1}{2\lambda} ( e^{\lambda  s^2} -1)\quad &\mbox{if}\,\, |s|\leq k\,,\\\\
\vpla (k) (|s|-k)+\frac{1}{2\lambda} ( e^{\lambda  k^2} -1) \quad &\mbox{if}\,\, |s|> k\,.
\end{array}
\right.
$$
Note that 
$$
 \vpla (k)  |s|  - e^{\la k^2}  \left(k^2-\frac{1}{2\lambda}\right)   -\frac{1}{2\lambda}   \leq \Phi_{\lambda , k} (s) \leq \vpla (k)  |s| \,,
$$
so that we deduce 
\be\label{time}
\begin{array}{c}
\dys \int_0^t
\langle (u_n)_t , \vpla ( T_k (u_n)) \xi \rangle\\
\dys \geq 
\vpla (k) \int_{B_n } |u_n (t,x)|    \xi 
 -    \left[ e^{\la k^2}  \left(k^2-\frac{1}{2\lambda}\right) +\frac{1}{2\lambda}\right]
  \mis\{B_{R+\rho}\} \\
  \dys -  \vpla (k) \int_{B_n }   | u_0 (x)| \xi  \,.
\end{array}
\ee
{
We notice that assumptions \rife{g3}--\rife{g2} imply the following growth condition  on $h$:\begin{equation}\label{natural}
\exists c_1 >0 \ \mbox{such that }\ h(\tau)\leq c_1 (\tau^{\frac{p}{p-1}}+1)\,,\forall \tau\in \re^+\,.
\end{equation}}
Moreover, since  $k>L$,  by \rife{g1}--\rife{g4}, we have
$$
\begin{array}{c}
\intbnp  g(t,x,u_n,\D u_n) \vpla  (T_k (u_n)) \xi\\ \dys 
= 
\int_{Q_n^t \cap\{ |u_n|\geq L\}}  g(t,x,u_n,\D u_n) \vpla  (T_k (u_n)) \xi 
\end{array}
$$
$$
\begin{array}{c}
\dys +
\int_{Q_n^t \cap\{ |u_n|\leq L\}}  g(t,x,u_n,\D u_n) \vpla  (T_k (u_n)) \xi
\geq 
\intbnp   h(|\D u_n|^{p-1}) |\vpla ( T_k (u_n)) |   \xi\\ \dys
-  \intbnp  \Big(\widetilde{\gamma}_k |\D T_k (u_n) |^p + |\widetilde{g_k} (t,x)| \Big)   | \vpla ( T_k (u_n)) | \xi\,,
\end{array}
$$
where $\widetilde{\gamma}_k ={\gamma_k}+c_1$ and $\widetilde{g}_k ={g_k}+c_1$, and $c_1$ is the constant appearing in \rife{natural}.
On the other hand by \rife{a2} and since $\xi$ can be chosen such that  $\xi= \sigma (\eta)$, where $\eta$ satisfies \rife{xi}, too,  and $\sigma$ is the function defined in Lemma \ref{cutoff}, we can apply  Proposition \ref{prop}  with $\delta=\frac{1}{2 \beta}$.
Thus there exists a constant $C=C (\lambda, k, \beta,T)$ such that  
$$
\begin{array}{c}
\intbnp a (t, x ,u_n, \D u_n )\cdot  \D \xi   \vpla ( T_k (u_n))   \\
\geq -\beta \intbnp  |\D u_n |^{p-1}  |\D \xi| |  \vpla ( T_k (u_n))  | \\
\geq - \intbnp  \frac12  h\left(  |\D u_n |^{p-1}\right)  |  \vpla ( T_k (u_n))  | \xi   - C  \mis \{ B_{R+\rho}\}\,.
\end{array}
$$
By substituting  the above inequalities into \rife{est1}, we deduce
\be \label{est11}
\begin{array}{c}
\dys \vpla(k) \int_{B_n }  |u_n (t,x)|   \xi 
 + \intbnp  | \D  T_k (u_n) |^p \Big[\alpha \vpla' ( T_k (u_n))
- \widetilde{\gamma}_k    |\vpla (T_k (u_n))| \Big] \xi \\[2.0 ex]
\leq 
  \mis\{B_{R+\rho}\} 
+ \dys \vpla(k) \left[ \int_{B_n }   | u_0 (x)| \xi 
+ \intbnp  \Big(|f_n(t,x)| + | \widetilde{g}_k (t,x)|\Big)   \xi\right] \,.
\end{array}
\ee

Note that  both $f_n (t,x)\xi$ and $\widetilde{g}_k (t,x)\xi $ are bounded in $L^1 (Q^T_R)$, therefore,  we choose $\lambda>\frac{\widetilde{\gamma}^2_k}{8\alpha^2}$ so that \rife{vpla} holds and we deduce  that there exists a constant (depending on $k$) such that 
\be\label{tronc}
\dys \sup_{t\in (0,T) }\int_{B_R }  |u_n (t,x)| +
\dys \int_{ Q_R^T } |\D T_k (u_n)|^p \leq C( k) \,,\qquad \forall k>0\,.
\ee

This implies, since obviously $\|T_{k}(u_{n})\|_{L^{p}((0,T)\times B_{R})}\leq C (R,T) k$,  that  $T_k (u_n) $ is bounded in $ L^p (0,T; W^{1,p} (B_R ))$, $\forall R>0$. Thus, up to subsequences (not relabeled) $T_k (u_n)$ weakly converges toward a function $v_k$ in $L^p (0,T; W^{1,p} (B_R ))$. Moreover the sequence $\{u_n\}$ is bounded in 
$L^{\infty} (0,T; L^{1} (B_R ))$\,.\\[2.0 ex]
Hence,  from \rife{tronc} we deduce (integrating between $0$ and $T$), $\forall j>0$, 
$$
\begin{array}{c}
j \mis \{(t, x) \in Q^T_R  \,:\, |u_n| \geq j \} 
\\
\dys \leq \int_{ \{(0,T)\times B_{R+\rho}\} \cap\{(t,x): |u_n| \geq j \} } |u_n (t,x) | \xi 
\leq \int_0^T \int_{B_n } |u_n (t,x) | \xi \leq CT\,,
\end{array}
$$
so that 
\be\label{misura}
\mis \{(t, x) \in Q_R^T \,:\, |u_n| \geq j \}   \leq \frac{C T}{j}\,.
\ee

Moreover, choosing $S'_k (u_n)\xi $ as test function in \rife{weakn} ($S_k$ has
been defined in \rife{sgei}), we deduce   that $ \big(S_k (u_n)\xi\big)_t$ is bounded in $L^1 (Q^T_R) + L^{p'} (0,T; W^{-1,p'} ( B_{R+\rho}))$ and so, using Corollary 4 in \cite{si}, we have that $S_k (u_n)\xi$ is strongly compact in $L^1 ((0,T)\times B_{R+\rho})$. Hence, up to subsequences (not relabeled), it converges a.e. as $n$ diverges. Using a diagonal argument, it follows that $u_n \to u$ for a.e. $(t,x)\in Q^T_R$, $\forall R>0$, and consequently there exists a measurable function $u (t,x)$ such that $u_n\to u$ a.e. in $(0,T)\times \rn$. Finally,  we note that \rife{misura} and the a.e. convergence of $u_n$ imply,  by Vitali's theorem,   that $u_n$ is compact in $L^1 (0,T; L^1_{\rm loc} (\rn))$ and consequently that  
$$
T_k (u_n)\rightharpoonup T_k (u) \quad \mbox{in}\,\, \eneloc\,,
$$
so that $v_k=T_k (u)$.

\medskip

\noindent{\bf Estimates on the lower order term.}
Let us choose, $\forall \eps>0$, 
$\psi =\frac{T_{\eps} ( u_n)}{\eps} \xi$ as test function in \rife{weakn}, so that we have:
$$
\dys \intbn \frac{d }{dt } \left( \frac{\Theta_{\eps}  (u_n)}{\eps} \right) \xi
+ \frac{1}{\eps} \intbn a (t, x ,u_n, \D u_n )\cdot  \D T_{\eps}(  u_n)    \xi 
$$ $$
+ \intbn a (t, x ,u_n, \D u_n )\cdot  \D  \xi  \frac{T_{\eps}  (u_n) }{\eps}  
$$
$$
+\intbn g(t,x,u_n,\D u_n) \frac{T_{\eps} ( u_n)}{\eps}   \xi = \intbn f_n (t,x) \frac{T_{\eps} ( u_n)}{\eps}\xi 
\leq \intbn | f_n(t,x)|     \xi 
\,,
$$
where $\Theta_k (s)= \int_0^s T_k (\tau)d\tau$, $\forall k>0$. 
We first note that 
$$
0 \leq  \frac{\Theta_{\eps} (s)}{\eps} \leq |s|, \quad \forall s\in \re, 
 $$
and by \rife{a1}, we deduce 
$$
\frac{\alpha  }{\eps} \intbn | \D T_{\eps}(u_n)|^p    \xi 
+ \intbn a (t, x ,u_n, \D u_n )\cdot  \D  \xi  \frac{T_{\eps} (u_n)}{\eps}
$$
$$
+\intbn g(t,x,u_n,\D u_n) \frac{T_{\eps} (u_n)}{\eps}  \xi
\leq \| f_n(t,x)\|_{L^1(Q^T_{R+\rho})}+ 
\int_{B_n } |u_n^0 (x)|  \xi
 \,.
$$
Moreover, using \rife{g1}--\rife{g4}, we have (as above  $\widetilde{\ga}_L ={\ga}_L +c_1$ and $\widetilde{g}_L ={g}_L +c_1$)
$$
\intbn g(t,x,u_n,\D u_n) \frac{T_{\eps} (u_n)}{\eps} \xi
$$ $$
\geq \frac12 \intbn h(|\D u_n|^{p-1}) \xi 
+ \frac12 \intbn |g(t,x,u_n,\D u_n) | \left|\frac{T_{\eps}(u_n)}{\eps}  \right|\xi 
$$ $$
-\widetilde{\ga}_L \int_{Q_{T}^{n} \cap \{|u_n|\leq L\}} |\D T_L (u_n)|^p \xi
- \int_{Q_{T}^{n} \cap \{|u_n|\leq L\}} |\widetilde{g}_L (t,x) |\xi \,,
$$
As we have already noticed, we can choose $\xi$ such that  $\xi=\sigma (\eta)$ and so both  $\eta$ and $\xi$ satisfy  \rife{xi}. Thus, by  using Proposition \ref{prop} with $\delta=\frac{1}{4\beta}$, we have 
$$
\left|\intbn a (t, x ,u_n, \D u_n )\cdot  \D  \xi  \frac{T_{\eps}}{\eps} (u_n)\right|
\dys  $$
$$
\dys \leq \frac14 \intbn h \left(  | \D u_n |^{p-1}\right) \xi + C\  \mis (B_{R+\rho}) \,.
$$
Hence, dropping positive terms,  we have
$$
\frac14 \intbn h(|\D u_n |^{p-1}) \xi 
+ \frac12 \intbn |g(t,x,u_n,\D u_n) | \left|\frac{T_{\eps}(u_n)}{\eps}  \right|\xi 
$$ $$
\leq \widetilde{\ga}_L \int_{B_n \cap \{|u_n|\leq L\}} |\D T_L (u_n)|^p \xi
+ \int_{B_n \cap \{|u_n|\leq L\}} |g_L (t,x)| \xi
$$
$$
+ \int_{B_n } |u_n^0 (x)|  \xi
+ C \mis (B_{R+\rho}) 
\,,
$$
and, by \rife{tronc} and \rife{g3},  the right hand side of the previous inequality is  uniformly  bounded  (with respect to $n$). 
Thus, letting $\eps \to0$, Fatou Lemma yields
\be\label{bound}
\dys
\int_{Q_R^T} h(|\D u_n |^{p-1})  +  \int_{Q_R^T} |g(t,x,u_n,\D u_n) |  
\leq C_R\,.
\ee

\noindent{\bf Equiintegrability of the lower order term and uniform estimates on  stripes.}
Let us choose
$\psi =  \gamma_j ( u_n)   \xi$, $\forall j>L$, in \rife{weakn}  
where 
$\gamma_j (s)=T_1(G_j(s))$, and moreover we denote by $\Gamma_j (s)=\int_0^s \gamma_j (t)dt$;  we note that 
\be\label{gamma}
|G_{j+1} (s)|\leq  \Gamma_j (s)\leq  |G_{j} (s)|\,.
\ee
Thus we have:
$$
\begin{array}{c}
\dys \int_0^T  \langle  (u_n)_t   \,,\,\gamma_j ( u_n)  \xi \rangle
+ \intbnp  a (t, x ,u_n, \D u_n )\cdot  \D  \xi  \gamma_j ( u_n)   \\[2.0 ex]
+ \intbnp  a (t, x ,u_n, \D u_n )\cdot  \D u_n \gamma_j' ( u_n )    \xi
  \\[2.0 ex]
+\intbnp  g(t,x,u_n,\D u_n)    \gamma_j ( u_n ) \xi
=\intbnp  f_n(t,x)    \gamma_j ( u_n  ) \xi\,.
\end{array}
$$

Thus, since $j> L$, using that $|\gamma (s)|\leq  1$ and \rife {a1} we get 
$$
\dys \int_{B_n } \Gamma_j (| u_n (t,x ) |)\xi
+ \intbnp a (t, x ,u_n, \D u_n )\cdot  \D  \xi   \gamma_j ( u_n)  $$ $$
+ \int_{Q_n^t \cap \{ j \leq  |u_n|\leq  j+1\}} a (t, x ,u_n, \D u_n )\cdot  \D u_n    \xi 
$$ $$
+\frac12 \int_{Q_n^t \cap \{ |u_n|\geq j\}}  h(|\D u_n|^{p-1})    |\gamma_j ( u_n )| \xi 
+\frac12 \int_{Q_n^t \cap \{ |u_n|\geq j\}}  |g(t,x,u_n,\D u_n)|   | \gamma_j ( u_n  )| \xi 
$$ $$
\leq  \int_{Q_n^t \cap \{ |u_n|\geq j \}} |f_n(t,x)|  \xi   
+ \dys \int_{B_n} \Gamma_j (| u_n^0 (x)|)\xi\,.
$$
On the other hand by \rife{a2}, and choosing $\xi=\sigma(\eta)$ as above, we deduce by Proposition \ref{prop} applied with $\delta=\frac{1}{2\beta}$, 
$$
 \int_{Q_n^t \cap \{ |u_n|\geq j\}}  |a (t, x ,u_n, \D u_n )\cdot  \D \xi | | \gamma_j (u_n)| 
$$
$$
\leq  
\frac12  \int_{Q_n^t \cap \{ |u_n|\geq j\}}  h  \left( | \D u_n |^{p-1}\right)  | \gamma_j (u_n)| \xi
$$ $$
 + C  \mis \{ (t,x)\in   (0,T)\times B_{R+\rho} \,:\, |u_n|\geq j \}\,, 
 $$
and the last term tends to $0$ (uniformly with respect to $n$) as $j$ diverges by \rife{misura}.
Moreover, by using \rife{gamma}  we deduce, dropping the positive term,  
$$
\dys \int_{B_n } G_{  j+1  } (|  u_n (t,x) |) \xi
+   \int_{Q_n^t \cap \{ j \leq  |u_n|\leq  j+1\}} a (t, x , u_n , \D u_n)\cdot \D u_n     \xi 
$$ $$
+\frac12 \int_{Q_n^t \cap \{ |u_n|\geq j\}}  |g(t,x,u_n,\D u_n)|  | \gamma_j ( u_n) | \xi 
$$ $$
\leq  \int_{Q_n^t \cap \{ |u_n|\geq j \}} |f_n(t,x)|  \xi   
+ \dys \int_{B_n} G_{ j} (|u_n^0 (x)|  )\xi 
+ \eps (j)
$$
Since both $u^n_0 (x)\xi $ and $f_n(t,x)\xi $ are strongly compact in $L^1 (B_{R+\rho})$ and $L^1 (Q^T_{R+\rho})$ respectively, we obtain,  dropping positive terms, 
\be\label{uniform}
\begin{array}{ll}
\dys \liminf_{j\to \infty} \, \sup_{n\in \na} \,    
\dys \left[ 
\dys \int_{Q_n^T \cap \{ j \leq  |u_n|\leq  j+1\}} a (t, x , u_n , \D u_n)\cdot \D u_n     \xi   \right. \\[2.0 ex] 
\qquad \qquad \left.  \dys +  \int_{Q_n^T \cap \{ |u_n|\geq j+1 \}}  |g(t,x,u_n,\D u_n)|   \xi \right] =0\,.
\end{array}
\ee

Note that the above estimate, in fact, allows us to say,  using \rife{g2} and since  $h(s)$ is superlinear at infinity, that
\be\label{equi}
 \sup_{n\in \na} \int_{Q_R^T \cap \{ |u_n| \geq j \}}  |\D u_n |^{p-1}  =\eps (j )\,.
\ee

\noindent{\bf Strong convergence of truncations.}
Let $\vpla (s)$ be the function introduced in \rife{vpla}, where $\lambda >0$ will be fixed in the sequel. We set $T_k (u)_{\nu}= \eta_{\nu} (T_k (u), T_k (u_0)) $, where  $\eta_{\nu} (\cdot)$  has been defined in Lemma \ref{lg}. 

Let us choose $\psi= \vpla ( z_{n,\nu}) S_j  '(u_n) \xi  $ as test function  in \rife{weakn}, where $z_{n,\nu}=T_k (u_n)-T_k (u)_{ \nu}$, $k\geq L$, and  $S_j (s)$ is {as} in  \rife{sgei}. 
Thus we have  
\be\label{nabo}
\begin{array}{c}
\dys \int_0^T  \langle  S_j (u_n)_t  \,,\, \vpla (z_{n,\nu} ) \xi \rangle
\dys+ \intbn a(t, x,u_n, \D u_n )\cdot  \D \xi \, \vpla (z_{n,\nu})  S'_j (u_n)
\\
+ \intbn a(t, x ,u_n, \D u_n )\cdot  \D  \Big(T_k (u_n)-T_k (u)_{ \nu}\Big)  \, \vpla' (z_{n,\nu}) S_j '(u_n) \xi  
\\
+ \intbn a (t, x ,u_n, \D u_n )\cdot  \D u_n \, \vpla (z_{n,\nu}) S_j ''(u_n) \xi  
\\
+\intbn g(t,x,u_n,\D u_n)  \vpla ( z_{n,\nu} )S_j '(u_n) \xi \dys
=\intbn f_n(t,x)  \vpla (z_{n,\nu} )S'_j (u_n) \xi\,.
\end{array}
\ee
We first note that 
$$
 \intbn a (t, x ,u_n, \D u_n )\cdot  \D  \Big( T_k (u_n)-T_k (u)_{ \nu}\Big)  \, \vpla' ( z_{n,\nu}) S_j '(u_n) \xi  
$$
$$
= \int_{Q_n^T \cap \{|u_n |\leq k\}}  a (t, x ,u_n, \D u_n )\cdot  \D  \Big( T_k (u_n)-T_k (u)_{ \nu}\Big)  \, \vpla' ( z_{n,\nu}) S_j '(u_n) \xi  
$$
$$
- \int_{Q_n^T \cap \{|u_n |\geq k\}}  a (t, x ,u_n, \D u_n )\cdot  \D  T_k (u)_{ \nu}  \, \vpla' ( z_{n,\nu} ) S_j '(u_n) \xi  \,.
$$
Using \rife{tronc} and recalling that Supp($S_{j}') \subset [-j-1, j+1]$, there exists $\varsigma_{k,j}\in (L^{p'} (Q_R^T))^{N+1} $ such that, 
$$
\dys \lim_{n\to \infty} \int_{Q_n^T \cap \{|u_n |\geq k\}}  a (t, x ,u_n, \D u_n )\cdot  \D  T_k (u)_{ \nu}  \, \vpla' (z_{n,\nu}) S_j '(u_n) \xi  
$$
$$
\dys = \int_{Q_{R+\rho}^T \cap \{|u |\geq k\}}   \varsigma_{k,j} \cdot  \D  T_k (u)_{ \nu}  \, \vpla' (z_{\nu}) \xi  \,,
$$
and last integral tends to $0$ as $\nu$ diverges. In fact we have that  $T_k (u)_{\nu} \to T_k (u)$ strongly in $L^p (0,T; W^{1,p}_{\rm loc} (\rn))$, and consequently 
$|\D T_k (u)_{\nu} |\chi_{\{|u|\geq k\}}$ tends to zero strongly in $L^p (0,T\,;\, L^{p}_{\rm loc} (\rn))$. 
Thus 
$$
 \intbn a (t, x ,u_n, \D u_n )\cdot  \D  \Big( T_k (u_n)-T_k (u)_{ \nu}\Big)  \, \vpla' (z_{n,\nu}) S_j '(u_n) \xi  
$$
$$
= \int_{Q_n^T}  a (t, x ,u_n, \D T_k (u_n) )\cdot  \D  \Big( T_k (u_n)-T_k (u)_{ \nu}\Big)  \, \vpla' (z_{n,\nu}) S_j '(u_n) \xi  
+\eps (n,\nu)\,.
$$
On the other hand, since $k\geq L$, 
$$
\intbn g(t,x,u_n,\D u_n)  \vpla (z_{n,\nu} )S'_j (u_n) \xi
$$ 
 $$
 \geq  \int_{Q_n^T} h (|\D u_n|^{p-1}) |\vpla (z_{n,\nu})| S'_j (u_n) \xi
-  \int_{Q_n^T \cap \{|u_n |\leq k\}} 
| \widetilde{ g }_k (t,x)|
 |\vpla (z_{n,\nu} )| \xi
$$ 
$$
-  \frac{\widetilde{ \ga}_k}{\alpha} \int_{Q_n^T \cap \{|u_n |\leq k\}} 
  a (t, x ,u_n, \D  T_k (u_n) )\cdot \D (  T_k (u_n)-T_k (u)_{\nu})  S'_j (u_n)
| \vpla (z_{n,\nu} )| \xi
 $$  $$
 +  \frac{\ga_k}{\alpha} \intbn  a (t, x ,u_n, \D  T_k (u_n) )\cdot \D T_k (u)_{\nu}  S'_j (u_n)
 \vpla (z_{n,\nu} ) \xi\,,
 $$ 
where, as before, $\widetilde{ \ga}_k=  { \ga}_k+c_1$ and $\widetilde{ g }_k=g_k+c_1$.
Reasoning as before, 
$$
 \frac{\ga_k}{\alpha} \intbn  a (t, x ,u_n, \D  T_k (u_n) )\cdot \D T_k (u)_{\nu}   \vpla (z_{n,\nu} )  S'_j (u_n) \xi =\eps(n,\nu)\,,
$$
and  
$$
\intbn \Big[|f_n(t,x)| + |\widetilde{ g }_k (t,x)| \Big] \vpla (z_{n,\nu} ) S'_j (u_n) \xi= \eps (n,\nu)\,.
$$
Gathering the above informations  together, we deduce:
\be\label{nabo2}
\begin{array}{c}
\dys \int_0^T  \langle  S_j (u_n)_t  \,,\, \vpla (z_{n, \nu} ) \xi \rangle
+ \intbn a (t, x ,u_n, \D  u_n )\cdot  \D \xi \, \vpla (z_{n, \nu})  S'_j (u_n)
\\
+ \intbn a (t, x ,u_n, \D T_k (u_n) )\cdot  \D  \Big(T_k (u_n)-T_k (u)_{ \nu}\Big)  \, \vpla' (z_{n,\nu}) S_j '(u_n) \xi  
\\
\dys+ \intbn a (t, x ,u_n, \D u_n )\cdot  \D u_n \, \vpla (z_{n,\nu}) S_j ''(u_n) \xi 
\\
-   \frac{\widetilde{ \ga }_k}{\alpha}  \intbn
 a (t, x ,u_n, \D  T_k (u_n) )\cdot \D \Big(  T_k (u_n)-T_k (u)_{\nu} \Big) 
 \vpla (z_{n,\nu} )S'_j (u_n)  \xi
 \\
+ \dys \int_{Q_n^T} h (|\D u_n|^{p-1}) |\vpla (z_{n,\nu})| S'_j (u_n) \xi\leq  \eps (n,\nu) \,.
\end{array}
\ee

Moreover, using that $\xi=\sigma (\eta)$, $\eta$ chosen as in \rife{xi} and thanks to  Proposition \ref{prop},   the second integral in \rife{nabo2} is estimated as
$$
\left|\intbn a (t, x ,u_n, \D u_n )\cdot  \D \xi \, \vpla (z_{n,\nu} )  \right|
\leq \frac{1}{2} \intbn h ( |\D u_n|^{p-1} )|\vpla (z_{n,\nu} )| S'_j (u_n) \xi  + \eps (n,\nu).
$$

Thus by \rife{nabo2}, dropping positive terms, we have
\be\label{nabo3}
\begin{array}{c}
\dys \int_0^T  \langle S_j (u_n)_t  \,,\, \vpla (z_{n,\nu} ) \xi \rangle
\\[1.5 ex] \dys
+ \intbn a (t, x ,u_n, \D T_k (u_n) )\cdot  \D  z_{n,\nu} 
\big[
 \vpla' (z_{n,\nu}) 
-  \frac{\widetilde{ \ga }_k}{\alpha}   
 \vpla (z_{n,\nu}) \Big] S'_j (u_n) \xi
\\[1.5 ex] \dys
\dys + \intbn a (t, x ,u_n, \D u_n )\cdot  \D u_n \, \vpla (z_{n,\nu}) S_j ''(u_n) \xi
\leq  \eps (n,\nu) \,.
\end{array}
\ee

Noticing that, by definition of $T_k (u)_{\nu}$,
$$
-\intbn a (t, x ,u_n, \D  T_k (u)_{\nu} )\cdot \D z_{n,\nu}
\Big[  \vpla' (z_{n,\nu}) -  \frac{\widetilde{ \ga }_k}{\alpha} |\vpla (z_{n,\nu} )|\Big] S'_j (u_n)  \xi=\eps (n,\nu)\,,
$$
we   can add  this  quantity  in both sides of \rife{nabo3}. Moreover,   by  \rife{uniform}, we get
$$
\left| \intbn a (t, x ,u_n, \D u_n )\cdot  \D u_n \, \vpla (z_{n,\nu}) S_j ''(u_n) \xi\right|\leq \eps(j)\,.
$$
Finally, in order to get rid of the integral involving the time derivative of $S_{j}(u_{n})$,  we apply the following inequality, whose proof is postponed at the end of this Section.\\[2.0 ex]
{\bf Claim.} $\forall j\geq j_0$:
\be \label{tempo}
\int_0^t \langle S_j (u_n)_t \, , \, \vpla (z_{n,\nu} )  \xi \rangle \geq \eps (n,\nu) \,.
\ee
\\[2.0 ex]
Using \rife{tempo} in \rife{nabo3} we deduce that, for $j$ large enough,
$$
\begin{array}{c}
\intbn (a (t, x ,u_n, \D T_k (u_n) )-a (t, x ,u_n, \D T_k (u)_{\nu} )) \cdot  \D   z_{n,\nu}    \\
\times \, \left[\vpla' (z_{n,\nu}) 
- \frac{\widetilde{ \ga}_k}{\alpha}  \vpla' (z_{n,\nu}) \right]
S_j '(u_n) \xi  
\leq  \eps (n,\nu) +\eps (j)\,.
\end{array}
$$
By a suitable choice of $\lambda$ (according with \rife{vpla} applied with $a=1$ and $b= \frac{\widetilde{ \ga }_k}{\alpha}$) we deduce that 
 $$
\begin{array}{c}
\intbn (a (t, x ,u_n, \D T_k (u_n) )-a (t, x ,u_n, \D T_k (u)_{\nu} )) \cdot  \D  z_{n,\nu} 
S_j '(u_n) \xi  
\leq  \eps (n,\nu)+\eps (j) \,.
\end{array}
$$
Lemma 5 in \cite{bmp} yields 
\be\label{convtro}
T_k (u_n) \to T_k (u)  \quad \mbox{strongly in } L^p (0,T; W^{1,p} (B_R )) \,.
\ee

Note that the above convergence implies that, up to subsequences, $\D T_k (u_n)$ a.e. converges to $\D T_k (u)$, and, by a diagonal argument,  we conclude that  (again up to not relabeled subsequences)
\be\label{aeconv}
\D u_n \to \D u \quad \mbox{a.e.}\,.
\ee
Moreover combining \rife{convtro} with \rife{equi} and \rife{aeconv} we deduce,  using Vitali Theorem, that 
\be \label{stro}
|\D u_n|^{p-1} \to |\D u |^{p-1} \quad \mbox{strongly in } \, L^1 (  0,T ; L^1_{\rm loc} (\rn))\,,
\ee
and by    \rife{uniform}, \rife{convtro} and \rife{aeconv} we have 
\be\label{cpt}
g(x, u_n , \D u_n)\xi  \to g(x,u, \D u )\xi  \quad \mbox{strongly in } \, L^1  (  (0,T) \times \rn )\,.
\ee
\\[2.0 ex]
\noindent{\bf Passing to the limit.}
Let us choose  $\psi =\phi (t,x) S' (u_n)$ in \rife{weakn} where  $\phi$ is in $C_0^{1} ([0,T)\times\rn)$ and $S$ is as in Definition \ref{def}. We have:

\be\label{pass}
\begin{array}{c}
\dys \int_0^T  \langle   (u_n)_t   \,,\, \phi S' (u_n)\rangle
+ \intbn a (t, x ,u_n, \D u_n )\cdot  \D  \phi   S'(u_{n}) \\[2.0 ex]
+ \intbn a (t, x ,u_n, \D u_n )\cdot  \D  u_{n} S''(u_{n})\phi \\[2.0 ex]
 +\intbn g(t,x,u_n,\D u_n)    \phi
=\intbn f_n(t,x)     \phi\,.
\end{array}
\ee

We first note that there exists $R>0$ such that supp $\phi (x,t)\subset (0,T)\times B_R$, so  that, integrating by parts,  and recalling that $\phi(T,x)=0$ we get: 
$$
\lim_{n\to +\infty}
\dys \int_0^T  \langle  (u_n)_t   \,,\, \phi S' (u_n) \rangle
$$
$$
=
\lim_{n\to +\infty}
\dys \int_{\rn}  S (u_n )  \phi (T,x)
-\dys \int_{\rn}  S (u_n^0 )  \phi (0,x)
- \int_{Q^T} S (u_n) \, \phi_t (t,x)
$$
$$
=
-\dys \int_{\rn}  S (u_0 )  \phi (0,x)
- \int_{Q^T}  S (u) \, \phi_t (t,x)\,,
$$ 
where $Q^T=(0,T)\times \rn$.
Moreover, by \rife{a2}  and \rife{stro}, we have 
$$
\lim_{n\to +\infty}
 \intbn a (t, x ,u_n, \D u_n )\cdot  \D  \phi S' (u_n) 
= \intbn a (t, x ,u , \D u )\cdot  \D  \phi S' (u) \,,
$$
while by \rife{convtro} and \rife{cpt}, \rife{uniform} we deduce  that both 
$$
\lim_{n\to +\infty}
 \intbn a (t, x ,u_n, \D u_n )\cdot  \D u_n \, S'' (u_n)  \phi 
= \intbn a (t, x ,u , \D u )\cdot  \D u  \, S'' (u )  \phi \,,
$$
and 
$$
\lim_{n\to +\infty}
\intbn g(t,x,u_n,\D u_n) S' (u_n)  \phi
=\intbn g(t,x,u,\D u) S' (u)  \phi\,.
$$
 Finally, since $f_n\xi \to f\xi $ in  $L^1 (Q^T)$ we can  pass to the limit in the last integral in \rife{pass}. Consequently  $u(t,x)$ is a solution for \rife{main} in the sense of Definition \ref{def}.

\medskip

\noindent{\bf Continuity with values in $L^1_{loc}$.} We   prove that $u_n$ converges up to subsequences (not relabeled) toward $u$ in $C^0([0,T];L^1_{loc} (\rn))$ by  using a classical monotonicity argument. 

For any integers $n$ and $m$, let us consider the formulations of \rife{weakn} with indexes  $n$ and $m$, respectively. For any $R,\rho>0$,  let us multiply the formulation of both $u_n$ and $u_m$ by $T_1(u_n - u_m)\xi $, with $n,m >R+\rho$, where $\xi$ has been defined   in \rife{xi}. By subtracting the two resulting identities, and dropping the energy terms (by  \rife{a3}), we get: 
$$
\begin{array}{c}
\dys \sup_{t\in [0,T]} \int_{B_R}\Theta_1 (u_n-u_m)(t,x)\\[1.5 ex]
\leq \|a(t,x,u_n, \nabla u_n)-a(t,x, u_m ,\nabla u_n) \cdot \nabla \xi \|_{L^1 (Q^T_{R+\rho})}\dys 
\dys+\|u^0_{n}- u^0_{m}\|_{L^1(B_{R+\rho})}\\[1.5 ex]
+\|g(t,x,u_n,\nabla u_n)-g(t,x,u_m,\nabla u_m)\|_{L^1(Q^T_{R+\rho})}+\|f_n - f_m\|_{L^1 (Q^T_{R+\rho} )}
\,.
\end{array}
$$
Note that by \rife{stro} and \rife{cpt} and using that both $\{u_n^0\}_{n\in \na}$ and $\{f_n\}_{n\in \na}$ are strongly compact in $L^1 (Q^T_{R+\rho} )$ and $L^1(B_{R+\rho})$,  respectively, the right hand side above converges to $0$ as $n$ and $m$ diverge. Hence it  follows that,   up to subsequences (not relabeled),   $u_n$ is a Cauchy sequence in $C^0([0,T];L^1 (B_R))$.

\medskip

To complete the proof of Theorem \ref{esistenza} we need to prove that inequality \rife{tempo} holds: the proof  follows the outlines of  Lemma 3.2 in \cite{do}.

 \proof[of  \rife{tempo}]. 
We recall that, by previous estimates,  
$T_k (u_n) $ converges to $ T_k (u)$ weakly in $ L^p (0,T ; W^{1,p}_{\rm loc} (\rn))
$. 
Here we {exploit} an approximation argument by using Lemma \ref{lg}. 
We set,  for every $\sigma>0$, $u_{n,\sigma}= \eta_{\sigma} (u_n, u_n^0) $;  we know that  $u_{n,\sigma}\in L^{p}(0,T; W^{1,p}_{0}(B_{n})) $,   $(u_{n,\sigma})_{t} \in L^{p}(0,T; W^{1,p}_{0}(B_{n}))$, and  moreover, both
$$u_{n,\sigma}\longrightarrow u_{n}\ \ \text{in}\ \  L^{p}(0,T; W^{1,p}_{0}(B_{n}))\,,$$
and
$$(u_{n,\sigma})_{t} \longrightarrow (u_{n})_{t} \ \ \text{in}\ \   L^{p'}(0,T; W^{-1,p'}(B_{n}))+ L^1 (Q_n^T)\,,$$ 
with $u_{n,\sigma}(0,x)=u_{n}^{0}$.

 This approximation argument will allow us to consider derivatives with respect to $t$ of the composition between Lipschitz functions and $u_{n,\sigma}$. 
Thanks to these properties we have that 
\be\label{sigma}
\begin{array}{c}
\dys 
\int_0^T \langle S_j (u_n)_t \, , \, \vpla (z_{n,\nu} )\xi \rangle\\
\dys =\lim_{\sigma\to 0}\int_0^T \langle S_j (u_{n,\sigma})_t \, , \, \vpla (T_k (u_{n,\sigma})-T_k (u)_{\nu} )\xi \rangle
\end{array}
\ee 
Our aim is to prove that
$$
\int_0^T \langle S_j (u_{n,\sigma})_t \, , \, \vpla (T_k (u_{n,\sigma})-T_k (u)_{\nu} )\xi \rangle\geq \eps(n,\nu)\,.
$$
Note that, for any $j>k$,  we can write 
$$
\begin{array}{c}
S_j (u_{n,\sigma})=   T_k(u_{n, \sigma}) + G_k (S_j (u_{n, \sigma}))
\end{array}
$$
thus, if we define $\phi_\lambda (s)=\int_0^s \varphi_\lambda$, we have 
$$
\begin{array}{c}\dys
\int_0^T \langle S_j (u_{n, \sigma})_t \, , \, \vpla (T_k (u_{n, \sigma})-T_k (u)_{\nu} )\xi \rangle
\\[1.5 ex] \dys
=\int_{B_R}  \phi_{\la} (T_k (u_{n, \sigma})-T_k (u)_{\nu} )(T) \xi - \int_{B_R}  \phi_{\la} (T_k (u_{n, \sigma})-T_k (u)_{\nu} )(0)\xi
\\[1.5 ex] \dys
+\int_0^T \langle G_k (S_j (u_{n, \sigma}))_t \, , \, \vpla (T_k (u_{n, \sigma})-T_k (u)_{\nu} )\xi \rangle
\\[1.5 ex] \dys
+\int_{Q_R^T}  \nu (T_k (u)-T_k (u)_{\nu} ) \vpla (T_k (u_{n, \sigma})-T_k (u)_{\nu} )\xi \,,
\end{array}
$$
where we  used  that ${(T_k (u)_{\nu}})_t = \nu (T_k (u)-T_k (u)_{\nu})$. {Since} both   $\phi_{\la} (s)>0$, $\forall s\in \re$ and  $T_k (u_{n}^0)\xi  \to T_k (u_0) \xi $ $\ast$-weakly in $L^{\infty} (B_R)$,  we deduce that 
$$
\int_{B_R}  \phi_{\la} (T_k (u_{n, \sigma})-T_k (u)_{\nu} )(T)\xi - \int_{B_R}  \phi_{\la} (T_k (u_{n, \sigma})-T_k (u)_{\nu} )(0) \xi \geq  \eps(n, \nu) \,.
$$
Moreover 
$$
 \int_{Q_R^T}  \nu (T_k (u)-T_k (u)_{\nu} ) \vpla (T_k (u_{n, \sigma})-T_k (u)_{\nu} ) \xi
$$
$$
= \int_{Q_R^T}  \nu (T_k (u)-T_k (u)_{\nu} ) \vpla (T_k (u)-T_k (u)_{\nu} )\xi + \eps( \sigma, n) \geq  \eps( \sigma, n)\,,
$$
since  $s\cdot s e^{\la s^2 } \geq 0$. Finally,  we deal with the term
$$
\int_0^T \langle G_k (S_j (u_{n, \sigma}))_t \, , \, \vpla (T_k (u_{n, \sigma})-T_k (u)_{\nu} ) \xi\rangle\,.
$$
Integrating by parts we deduce that 
$$
\begin{array}{c}
\dys
\int_0^T \langle G_k (S_j (u_{n, \sigma}))_t \, , \, \vpla (T_k (u_{n, \sigma})-T_k (u)_{\nu} ) \xi \rangle
\\[1.5 ex] \dys
= \int_{B_R} G_k (S_j (u_{n, \sigma}))(T)  \vpla (T_k (u_{n, \sigma})-T_k (u)_{\nu} ) (T)\xi
\\[1.5 ex] \dys
- \int_{B_R} G_k (S_j (u_{n, \sigma}))(0)  \vpla (T_k (u_{n, \sigma})-T_k (u)_{\nu} ) (0) \xi
\\[1.5 ex] \dys
- \int_0^T \langle G_k (S_j (u_{n, \sigma})) \vpla' (T_k (u_{n, \sigma})-T_k (u)_{\nu} ) \xi \, ,\, (T_k (u_{n, \sigma})-T_k (u)_{\nu} )_t \rangle\,.
\end{array}
$$
Thus the first term on the right hand side is positive since
$$
\begin{array}{l}
\dys \int_{B_R} G_k (S_j (u_{n, \sigma}))(T)  \vpla (T_k (u_{n, \sigma})-T_k (u)_{\nu} ) (T) \xi
\\
=\dys \int_{B_R \cap \{u_{n, \sigma}>k\}} G_k (S_j (u_{n, \sigma}))(T)  \vpla (k-T_k (u)_{\nu} ) (T) \xi
\\
 + \dys\int_{B_R \cap \{u_{n, \sigma}<-k\}} G_k (S_j (u_{n, \sigma}))(T)  \vpla (-k-T_k (u)_{\nu} ) (T) \xi
\geq 0\, ,
\end{array}
$$
while the second term vanishes passing to the limit with respect to $\sigma$, $n$ and then $\nu$. 
Concerning  the last one, we note that 
since $G_k ( S_j (u_{n, \sigma}))$ is $0$ if $|u_{n, \sigma}|\leq k$, thus 
$$
- \int_0^T \langle G_k (S_j (u_{n, \sigma})) \vpla' (T_k (u_{n, \sigma})-T_k (u)_{\nu} )\xi \, ,\, (T_k (u_{n, \sigma})-T_k (u)_{\nu} )_t \rangle
$$
$$\dys 
= \nu \int_{Q_R^T} G_k (S_j (u_{n, \sigma})) \vpla' (T_k (u_{n, \sigma})-T_k (u)_{\nu} )(T_k (u)- T_k (u)_{\nu}  )  \xi\,.
$$
Finally, taking the limit respectively in $\sigma$ and  $n$,  we have 
$$
\begin{array}{c}
\dys  - \int_0^T \langle G_k (S_j (u_{n, \sigma})) \vpla' (T_k (u_{n, \sigma})-T_k (u)_{\nu} )\xi \, ,\, (T_k (u_{n, \sigma})-T_k (u)_{\nu} )_t \rangle
\\[1.5 ex]
\dys 
= \nu  \int_{Q_R^T} G_k (S_j (u)) \vpla ' (T_k (u)-T_k (u)_{\nu} )(T_k (u)- T_k (u)_{\nu}  ) \xi+\eps (\sigma,n)
\\[1.5 ex]
\dys 
= \nu\int_{Q_R^T\cap\{u>k\}} G_k (S_j (u)) \vpla ' (k-T_k (u)_{\nu} )(k- T_k (u)_{\nu}  ) \xi
\\[1.5 ex] \dys
+ \nu\int_{Q_R^T\cap\{u<-k\}} G_k (S_j (u)) \vpla ' (-k-T_k (u)_{\nu} )(-k - T_k (u)_{\nu}  )\xi  + \eps(\sigma,n)\geq \eps(\sigma,n)
\end{array}
$$
since $ \vpla ' (s)>0$, $\forall s\in \re$  and the {claim} is proved because of \rife{sigma}. 

\medskip 

This concludes the proof of Theorem \ref{esistenza}.
\qed
\medskip
\proof[of Theorem \ref{exls}]. 
 Let $u_n$ be a weak solutions of the following problem 
\begin{equation}\label{appls}
\left\{
\begin{array}{ll}
(u_n)_{t}- \dive a (t, x ,u_n,\D u_n)   +g(t,x,u_n,\D u_n) = f_n(t,x) \, &\mbox{in}\,Q^T_{\Omega}  ,\\[1.5 ex]
 u_n (t,x)= n & \mbox{on}\,   \partial_{\mathcal{P}} Q^T_{\Omega}  \,, \\[1.5 ex]
u_n (x,0)= u_n^0 (x) &\mbox{in}\,   \Omega ,
\end{array}\right.
\end{equation}
where $f_n (t,x)= T_n (f(t,x))$ and $ u_n^0 (x) =T_n (u_0 (x))$. The existence of  a weak solution for \rife{app} is still a consequence of the result of \cite{do}. This means that there exists a function $u_n$ such that  $u_n-n \in L^p (0,T ; W^{1, p}_0 (\Omega))$,  $(u_n)_t \in L^{p'} (0,T ; W^{-1, p' } (\Omega))$, $g(t,x,u_n,\D u_n) \in L^1 ((0,T)\times \Omega)$, and the following identity holds true 
\be\label{weakls}
\begin{array}{c}
\dys \int_0^T  \langle  (u_n)_t   \,,\, \psi \rangle
+ \int_{Q_{\Omega}^T} a (t, x ,u_n, \D u_n )\cdot  \D  \psi  \\[1.5 ex]
\dys +\int_{Q_{\Omega}^T} g(t,x,u_n,\D u_n)    \psi
=\int_{Q_{\Omega}^T} f_n    \psi\,, \\[1.5 ex]
\forall \psi \in L^p (0,T ; W^{1, p}_0 (\Omega))\cap L^{\infty} (Q_{\Omega}^T)\,.
\end{array}
\ee

\medskip
The idea of the proof is similar to the one of Theorem \ref{esistenza}. 
The main difference relies on the fact that now  we need to have an information about $u_n$ (and consequently $u$) at the boundary, and so we need first  to prove a \emph{global} (i.e. on the whole $\Omega$) estimate on the truncates in the energy space.
On the other hand, for the second part of the proof,   we follow  exactly the same outline of the one of Theorem \ref{esistenza}. Indeed, the estimates proved there are localized in $(0,T)\times B_R$, $\forall R>0$. Since, in order to pass to the limit in the equation, we need to use such estimates on any compact subset $\varpi \subset \subset (0,T)\times \Omega$, we observe that there exists $\omega \subset \subset \Omega$ such that 
$\varpi \subset \subset (0,T)\times \omega$. Thus $$
\exists M \in \na\,, \, x_i \in \Omega, \, r_i>0,\,  i=1,...,M,\,\mbox{  such that }
\omega \subset \bigcup_{i=1}^M B_{r_i} (x_i).
$$ 
It is clear that it is enough to prove all the estimates on a ball and without loss of generality we can suppose that  it is centered at the origin.

\medskip

\noindent{\bf Global estimate on truncations.} Let us choose,  $\forall n\geq k {\geq L} $,  $\psi = \vpla (T_k (u_n)-k)$ as test function  in \rife{weakls}, with $\lambda >0$ to be fixed later.  Thus we have 
$$
\begin{array}{c}
\intO \Upsilon_{\lambda, k} (u_n (t,x)) - \intO \Upsilon_{\lambda, k} (u_n^0 (x) ) \\
\dys + \int_{Q_{\Omega}^T} a (t, x ,u_n, \D u_n )\cdot  \D T_k (u_n) \vpla' (T_k (u_n)-k) \\
\dys +\int_{Q_{\Omega}^T} g(t,x,u_n,\D u_n)    \vpla (T_k (u_n)-k)\\
\dys =\int_{Q_{\Omega}^T} f_n^+    \vpla (T_k (u_n)-k) - \int_{Q_{\Omega}^T} f_n^-    \vpla (T_k (u_n)-k)\,,
\end{array}
$$
where 
$$
\Upsilon_{\lambda, k} (s) =
 \left\{
 \begin{array}{ll}
 -2k   e^{4\lambda k^2} (s+k) + \frac{1}{ 2\lambda} \left[  e^{4\lambda k^2} -  e^{\lambda k^2}\right] \quad &\mbox{if}\, s< -k\,,\\
 \frac{1}{ 2\lambda} \left[  e^{\lambda (s-k)^2} -  e^{\lambda k^2}\right] &\mbox{if}\, -k\leq s< k\,,\\
\frac{1}{ 2\lambda} \left[  1 -  e^{\lambda k^2}\right] &\mbox{if}\,  s\geq k\,,
 \end{array}
 \right. 
$$\
is a primitive of $\vpla (T_k (s)-k)$.
Note that, since $\Upsilon_{\lambda, k} (s)$ is decreasing and $\Upsilon_{\lambda, k} (0)=0$, then 
$$
\begin{array}{c}
\intO \Upsilon_{\lambda, k} (u_n (t,x)) - \intO \Upsilon_{\lambda, k} (u_n^0 (x)) 
\\ \dys
\geq \dys \int_{\Omega \cap \{0\leq u_n\leq k\}} \Upsilon_{\lambda, k} (u_n (t,x)) + 
 \dys \int_{\Omega \cap \{ u_n> k\}} \Upsilon_{\lambda, k} (u_n (t,x))
\\ \dys
  -  \dys \int_{\Omega \cap \{ u^0_n\leq -k\}} \Upsilon_{\lambda, k} (u_n^0 (x)) 
  -  \dys \int_{\Omega \cap \{ -k\leq u^0_n\leq 0 \}} \Upsilon_{\lambda, k} (u_n^0 (x)) 
\\[2.0 ex] \dys
\geq 
-\Big( \frac{1}{ \lambda} \left[    e^{\lambda k^2}-1\right]  
+ \frac{1}{ \lambda} \left[  e^{4\lambda k^2} -  e^{\lambda k^2}\right]
\Big) |\Omega| 
-2k   e^{4\lambda k^2} \|(u_n^0)^- \|_{\elle1}\,.
\end{array}
$$

Thus, by \rife{a1}, \rife{g31} and the assumptions on $f$ we deduce, since the function $\vpla (T_k (s)-k) \leq 0 $, $\forall s \in \re$, 
$$
\begin{array}{c}
\dys \alpha \int_{Q_{\Omega}^T} | \D T_k (u_n) |^p \vpla' (T_k (u_n)-k) 
\dys - \int_{Q_{\Omega}^T} \gamma_k |\D T_k (u_n)|^p    |\vpla (T_k (u_n)-k)|\\
\dys \leq   \vpla (2k) \int_{Q_{\Omega}^T} \Big[|f_n^- | +|g_k (t,x)|\Big]
+2k   e^{4\lambda k^2} \|(u_n^0)^- \|_{\elle1}\\
+
\Big( 
\frac{1}{ 2\lambda} \left[    e^{\lambda k^2}-1\right]  
+ \frac{1}{ \lambda} \left[  e^{4\lambda k^2} -  e^{\lambda k^2}\right]
+ 2k^2   e^{4\lambda k^2}
\Big) |\Omega|  \,.
\end{array}
$$
By fixing a suitable $\la>0$, so that  \rife{vpla} holds for $\vpla(s)$, we deduce that $k-T_k (u_n)$ is bounded in $\ene$ and so, up to subsequences (not relabeled),  it weakly converges in $\ene$. 

\medskip
As already pointed out, the conclusion of the Theorem follows exactly using the same steps of  Theorem \ref{esistenza}. \qed

\setcounter{equation}{0}
\section{Further Regularity}

In this section we are going to describe some local regularity properties for a renormalized solution of problem
\begin{equation}\label{mains}
\left\{
\begin{array}{ll}
u_t - \dive a (t, x ,u,\D u)   +g(t,x,u,\D u) = f(t,x) \quad &\mbox{in}\, (0,T) \times \Omega \\[1.5 ex]
u(0,x)= u_0 (x) &\mbox{in}\,   \Omega \,,
\end{array}\right.
\end{equation}
where $\Omega$ is a, possibly unbounded, domain in $\rn$. 

Let us first emphasize that in this section we would like to be able to choose test functions of the type $S'(u)\psi$ with $S'$ not compactly supported on $\re$ and such that $\psi(T,x)\neq 0$. 
In principle, according to Definition \ref{def}, we are not allowed to do that.  Anyway, after suitably modifying our definition, this fact can be made rigourous by an easy density argument.  
In fact, we can  choose $S'(u)= S'_j (u)M (u)$  where $M$ is a Lipschitz function and $S_j$ is defined in \rife{sgei}, in the renormalized formulation. 
Then,  we take the limit as $j$ diverges and we observe that  $S'_j (u)$ converges to $1$ both a.e. and $\ast$-weak in $L^{\infty}(Q_\Omega^T)$, and the term involving $S''$ vanishes thanks to \rife{fettine}.

Moreover, we need to deal with cut-off functions which do not depend on time; to do that we  choose a family  of functions of the type $\phi_{\delta} (t,x)=\xi(x)\psi_\delta (t)$ such that  $\xi(x)\in C^1_0 (\Omega)$ and $\psi_\delta(t)\in  C^1_0 ([0,T))$ that converges to $\chi_{[0,\tau]}$. 
Thus, according  Proposition \ref{marc}, a standard choice of $\psi_\delta(t) $ allows us to deduce   the following equivalent formulation that is the useful one in order to  {obtain} our regularity estimates: 
\be \label{renormr}
\begin{array}{c}
\dys 
\int_{\Omega} \mathcal{M} (u (\tau,x)) \xi (x)
+ \int_{Q_\Omega^{\tau}}  a (t, x ,u, \D u )\cdot \D u  M' (u) \xi   \\[2.0 ex]
\dys + \int_{Q_\Omega^{\tau}}  a (t, x ,u, \D u )\cdot   \D \xi M (u)  \\[2.0 ex]
\dys +\int_{Q_\Omega^{\tau}}  g(t,x,u ,\D u ) M (u) \xi
=\int_{Q_\Omega^{\tau}}  f (t,x) M (u) \xi  
+\int_{\Omega} \mathcal{M}  (u_0) \xi (x)
\,,
\end{array}
\ee
for any $0<\tau \leq T$, $\xi(x)\in C^1_0 (\Omega)$, and $\mathcal{M}' (s)=M (s)$, with  $\mathcal{M} (0)=0$.
\medskip

Finally,  we observe that, since the estimates we are going to prove in this section are localized, we will proceed as follows. We fix a ball (without loss of generality,  centered at 0) of radius $R$ contained in $\Omega$. Thus there exists $\rho>0$ such that $B_{R+\rho} \subset \subset \Omega$ and we will prove the estimates in $Q_R^T$, depending on quantities computed on $Q^T_{R+\rho}$. By covering any compact $\omega \subset \subset \Omega$ with balls we then obtain the results.
\medskip 

We start proving Proposition \ref{marc}.

\proof[of  Proposition \ref{marc}]. 
According to the formulation \rife{renormr} we are allowed to choose $\psi (t,x)=\xi (x)$, where $\xi$ is chosen as in \rife{xi} and such that Proposition \ref{prop} holds true,  and $M(s)=T_k (s)$, $\forall k\geq L$. Thus we have, recalling that $\Theta_k (s)=\int_0^s T_k (\tau)d\tau$
$$
 \int_{B_{R+\rho}} \Theta_k (u(x,t)) \xi^p  +
\frac{\alpha}{2^{p-1}} \int_{Q_{R+\rho}^T} |\D (T_k (u) \xi) |^p  +
 \frac12 \int_{Q_{R+\rho}^T} h (|\D u|^{p-1}) T_k (u)\xi^p  
 $$
$$
\leq k \|f\|_{L^1 (Q_{R+\rho}^T)}  + C_0 + \alpha \int_{Q_{R+\rho}^T} 
|T_k (u)|^{p} |\D \xi|^p + k \|u_0 \|_{L^1(B_{R+\rho})}\,.
$$
Note that, since Proposition \ref{p-1} holds true, then $u^{p-1}$ belongs to $L^1(Q_R^T)$ and so the last integral can be estimated by  $Ck$, for suitable $C>0$. Thus we deduce, by dropping positive terms,  
$$
 \int_{Q_R^T} |\D T_k (u)  |^p \leq C (k+1)\,,
$$
and so we conclude applying Lemma \ref{gajardo}.  
 Moreover the embedding of Marcinkiewicz spaces into Lebesgue ones (see \rife{embe'}) together with the assumption $p>2-\frac{1}{N+1}$ allow us to say that $u$ is bounded in $L^r (0,T;W^{1,r}_{\rm loc} (\Omega))$ for some $r>1$,  while from the equation in \rife {mains} (which is satisfied in the sense of distributions) we observe have $u_t\in L^{1}_{loc}(Q^T)+ L^{r'} (0,T;W^{-1,r'}_{loc} (\Omega))$. 
Thus the continuity with values in $L^{1}_{\rm loc}$ is an easy consequence of Theorem 1.1 of \cite{via}.
\qed

\medskip

Finally we give an idea of the proof of Theorem \ref{esse}.

\medskip

{\it Sketch of the proof  of Theorem \ref{esse}.} 
We start by giving the proof of   $(i)$ and $(ii)$. 

\medskip

 Let us fix any $0<R<R+\rho$ and consider  $B_R\subset B_{R+\rho} \subset\subset \Omega$.  Let us choose $M(s)= v_{\vare,j} (s)$, and $\psi=\xi^{\lambda}$ in \rife{renormr}, where $\lambda=\max \{p, \frac{q'p}{q'p -1},  \frac{q'p'}{q'p' -1}\}$,  $\xi (x)$ is as in \rife{xi} and 
$$
v_{\vare,j} (s) = \left[(|T_j(s)|+\eps)^{\gamma}-\vare^{\gamma}\right] \sign s, 
$$
 for any $0<\gamma \leq \overline{\gamma}$, with
\begin{equation}\label{gamma0}
\dys 
\overline{\gamma}= 
\left\{
\begin{array}{ll}
\frac{Nm(q-1)+q (m-1)[p(N+1)-2N]}{Nm -pq(m-1)}  \quad & \mbox{if}\,\, \rife{327} \, \, \mbox{holds}\\[2.0 ex]
\frac{N(p-1)(q-1)}{N-pq} & \mbox{if}\,\, \rife{327n} \, \, \mbox{holds}\,.
\end{array}\right.
\end{equation} 

We follow the same ideas of previous estimates, using the ellipticity condition, assumption \rife{g2} and Proposition \ref{prop}, and  finally letting $\eps$ tend to zero. Thus, we deduce that  there exists  a constant $C=C(\alpha,\beta,L, N,R,\rho, m,q)>0$,  but independent  of $j$, such that, 
$$
\begin{array}{c}
\dys 
\| \xi^{\frac{\lambda (\gamma+p-1)}{p(\gamma+1)}}  T_j (u) \|_{L^{\infty} (0,T; L^{\gamma+1} (B_{R+\rho}))}^{\gamma+1}
+
\int_0^T \| \xi^{\frac{\lambda }{p}}  |T_j (u)|^{\frac{\gamma+p-1}{p}} \|_{L^{p^{\ast}} (B_{R+\rho})}^p\\[2.0 ex] \dys 
\leq 
C\left[
\| f \|_{L^m (0,T; L^{q} (B_{R+\rho}))} \left(\int_0^T \left( \int_{B_{R+\rho}} \xi^{\frac{\lambda}{p}} |T_j (u)|^{\gamma q'}\right)^{\frac{m'}{q'}}\right)^{\frac1{m'}}\right.
\\[2.0 ex] \dys 
\left.
\dys +\int_{Q_{R+\rho}^T}  \xi^{\lambda-p}  |T_j (u)|^{p+\gamma-1}
 +\int_{B_{R+\rho} }    | u_0|^{ \gamma+1} +1
\right]\,,
\end{array}
$$
where we have applied a space-time H\"{o}lder inequality on the term involving the datum $f(t,x)$. 
Using the interpolation
 inequality with respect to the space variable and Young inequality with respect to the time variable, we deduce both that 
\be\label{nab-n}
\begin{array}{c}
\dys 
\|  T_j (u) \|_{L^{\infty} (0,T; L^{\gamma+1} (B_R))}
\leq 
C_1 \left[1+
\dys  \|T_j (u)\|^{p+\gamma-1}_{L^{p+\gamma-1} (Q_{R+\rho}^T)}
\right]^{\frac{1}{\gamma+1}}\,,
\end{array}
\ee
and 
\be\label{nab-n+1}
\begin{array}{c}
\|  T_j (u)\|_{L^{\gamma+p-1} (0,T;L^{{\frac{\gamma+p-1}{p}} p^{\ast}} (B_R))}\leq 
C_1 \left[1+
\dys  \|T_j (u)\|^{p+\gamma-1}_{L^{p+\gamma-1} (Q_{R+\rho}^T)}
\right]^{\frac{1}{\gamma +p -1 }}\,,
\end{array}
\ee
where $C_1$ also depends on $\| f \|_{L^m (0,T; L^{q} (B_{R+\rho}))}$ and, since $\gamma \leq \overline{\gamma}$, on  $\|u_0\|_{L^{\ol{\gamma}+1} (B_{R+\rho})}^{\ol{\gamma}+1}$,  too. 
Note that $\ol{\gamma}+1$ is the best summability we can expect as far as   the initial datum is concerned. 
By \rife{nab-n}, \rife{nab-n+1} and by applying  inequality \rife{interpolo} to the function $|T_j(u)|^{\frac{\gamma+p-1}{p}}$, we have:
\be\label{interpol}
\begin{array}{c}
\dys 
\|  T_j (u)  \|_{L^{\frac{p+\gamma-1}{p}\eta} (Q_R^T)}^{\frac{p+\gamma-1}{p}}
\\ \leq C
\|  T_j (u) \|^{1-\theta}_{L^{\infty} (0,T; L^{\gamma+1} (B_R))}\|  T_j (u) \|^{\theta}_{L^{\gamma+p-1} (0,T;L^{{\frac{\gamma+p-1}{p}} p^{\ast}} (B_R))}\,,
\end{array}
\ee
where $\eta$ and  $\theta$ satisfy 
$$
\frac{1}{\eta}=\frac{\theta}{p^\ast}+\frac{(1-\theta)({\gamma+p -1})}{{p(\gamma+1)}}, \ \ \ \frac{p}{\theta}\geq \eta. 
$$
It is easy to see that, if  
$$
\eta=p \left[1+\frac{p(\gamma+1)}{N(\gamma+p-1)}\right]
$$ then the above constraints are optimized. 
Thus gathering together \rife{interpol}, \rife{nab-n} and \rife{nab-n+1}
we deduce that there exists $C>0$ such that 
\be\label{interpol3}
\begin{array}{c}
\dys 
\|  T_j(u)  \|_{L^\frac{\eta(\gamma+p-1)}{p} (Q_R^T)}\leq 
C \left[1+
\dys  \|T_j(u)\|^{\gamma+p-1}_{L^{\gamma+p-1} (Q_{R+\rho}^T)}
\right]^{\frac{p}{\gamma+p-1} \max \{\frac{1}{\gamma+p-1} ,\frac{1}{\gamma+1}\}}\,.
\end{array}
\ee
Therefore, we control the norm  of $T_j (u)$ in $L^\frac{\eta(\gamma+p-1)}{p}$ of a cylinder  with the norm in $L^{\gamma+p-1}$ of a slightly larger cylinder. Moreover  such estimate is uniform with respect to $j$. Noticing   that $\eta>p$, in order to  conclude it is enough to perform an iteration method. We can construct  both $\overline{k}+1$ radii $0=\rho_0<\rho_1, ..., \rho_{\overline{k} -1}<\rho_{\overline{k}} =\rho$ and $\overline{k}+1$ exponents $\gamma_0<\gamma_1, ..., \gamma_{\overline{k}-1}<\gamma_{\overline{k}}=\overline{\gamma}$, such that  
$$
\gamma_0 +p-1<p-1+\frac{p}{N}\,,
$$
and $\frac{\eta(\gamma_{\overline{k}}+p-1)}{p}$ is our desired summability. 

Thus, applying \rife{interpol3} $\overline{k}+1$ times and  using Proposition \ref{marc},  we get the result.    To  deal with  the different time-space summability  stated in Theorem \ref{esse} we can argue in a similar way by applying  H\"{o}lder inequality. 

\medskip 

\noindent Now we deal with part $(iii)$ of the Theorem.
Let us denote by
\begin{equation}\label{notation}
A_{k,r} =\{x\in B_{\rho}(x_0):\ |u(t,x)|>k\},\ \forall r>0,
\end{equation}
and let us set, for any  fixed $\delta \in (0,1)$,
\be\label{lambda}
\begin{array}{c}
\dys  t_1 = \left[ \frac{ 1-\delta}{|B_{R+\rho}|^{\lambda_1}}\right]^{\lambda_2} \\[1.5 ex]
\dys  \lambda_1= \left(1-\frac{\sigma}{q(\sigma-p)}\right) \left(1-\frac{p}{\sigma}\right)
 \quad \mbox{and}\quad
\lambda_2= \frac{m \mu }{m (\mu-p)-\mu}\,,
\end{array}
\ee
where
\be\label{musigma}
\sigma=p\frac{Nm' + pq'}{Nm'},\ \ \ \ \mu=p\frac{Nm' + pq'}{Nq'}\ .
\ee

Let us choose $M(u)=T_j (G_{k} (u))$ in \rife{renormr},  $j>k>L$ and $\psi = \xi^p$ ($\xi $ chosen as in \rife{xi}) in the cylinder of \emph{height}  $t_{1}$,  where we will fix $\delta$ (and so $t_1$) later.

We also  recall that for $\xi$ chosen as in \rife{xi}, we have $|\D \xi|\leq \frac{c}{\rho}$. Thus, by standard computations,
we have
\begin{equation}\label{steppe1}
\begin{array}{l}
  \dys \iakrt |\nabla {T}_j (G_{k} (u))|^{p}  \xi^p
  \leq
\dys \frac{c_1}{\rho} \int_{ Q_R^{t_1} } | {T}_j (G_{k} (u))|\xi^{p-1}\\
\dys
+    \iakrt |f| | {T}_j (G_{k} (u))|^p\xi^p
+ \iakrt |f|\xi^p
\,.
\end{array}
\end{equation}

Moreover, we choose   $M_\eps (u)=[(|T_j ( {G}_{k}  (u))|+\eps)^{ p-1}-\eps^{p-1}] \sign u$,  $j>k>L$ and $\psi=\xi^p$ in \rife{renormr}. Thus, dropping positive terms, as  $\eps$ goes to zero, we get
\begin{equation}\label{steppe2}
\begin{array}{c}
  \dys \sup_{t\in (0,t_1)}\iakr \Theta_j (|G_k(u)|^p) (t) \xi^p
\leq
  \dys \frac{c_2}{\rho} \iakrt |T_j ( {G}_{k}  (u))|^{ p-1}  \xi^{p-1}
  \\
  +
\dys
c_3  \iakrt |f|  \xi^p
+\dys c_4  \iakrt |f|   |T_j ( {G}_{k}  (u))|^p  \xi^p\,,
\end{array}
\end{equation}
where, as before, $\Theta_j (s)$ denotes the primitive of $T_j (s)$ such that $\Theta_j (0)=0$.
Now we  apply Corollary \ref{lemma310}  with $w=T_j (|G_k (u)|) \xi$,   $\Omega=B_{R+\rho}$ and $T=t_1$. Thus,
for all $\mu$ and $\sigma$ satisfying \rife{330},    we deduce, by adding  \rife{steppe1} and \rife{steppe2}, and by applying Corollary \ref{lemma310}
\begin{equation}\label{firsti}
\begin{array}{c}
  \dys\left[\int_{0}^{t_{1}}\left(\int_{A_{k,R+\rho}} (|T_j (\gku)| \xi )^{\sigma}\right)^{\frac{\mu}{\sigma}}\right]^{\frac{p}{\mu}}
\leq
 \frac{c_5}{\rho^p }\iakrt  |T_j (\gku)|^{p_0}\\
\dys+ c_6  \iakrt |f|   |T_j ( {G}_{k}  (u))|^p  \xi^p
\dys +c_7  \iakrt |f|  \xi^p
\, ,
\end{array}
\end{equation}
with $p_0=\max \{ 1, p-1\}$.
Recalling the definitions of $\mu$ and $\sigma$ (see \rife{musigma}) and noticing that   both of them are greater than $p$,  we use H\"{o}lder inequality to estimate   the right hand side of \rife{firsti}, so that
$$
\begin{array}{c}
\dys  \iakrt |f|  |T_j ( {G}_{k}  (u))|^{p}  \xi^p \\
\dys \leq   |B_{R+\rho}|^{\lambda_1} t_1^{\frac{1}{\lambda_2}}
\|\xi  T_j (u)\|_{L^{\sigma} (0,t_1; L^{\mu} (A_{k,R+\rho}))}^p
\|  f \|_{L^{m} (0,t_1; L^{q} (B_{R+\rho}))}\,,
\end{array}
$$
where $\lambda_1$ and $\lambda_2$ have been defined in $\rife{lambda}$. We fix now $\delta$ (and consequently  we fix $t_1$) such that $c_6\|f\| (1-\delta)<\frac12$; note that this quantity depends on the data of the problem but not on $u$. Thus   from \rife{firsti} we deduce
\be\label{sarca}
\begin{array}{c}
  \dys
  \left[\int_{0}^{t_{1}}\left(\int_{A_{k,R+\rho}} |T_j (\gku)|^{\sigma} \xi^{\sigma} \right)^{\frac{\mu}{\sigma}}\right]^{\frac{p}{\mu}}\\
  \dys  \leq    \iakrt | f|\xi^p + \frac{c}{\rho^p }\iakrt  |T_j (\gku)|^{p_0}\,.
\end{array}
\ee
Moreover,  by H\"{o}lder  inequality it follows   that, for every $h>k$,
$$
\left[\int_{0}^{t_{1}}\left(\iakr|T_j (\gku)|^{\sigma} \xi^{\sigma} \right)^{\frac{\mu}{\sigma}}\right]^{\frac{p}{\mu}}
\geq (h-k)^{  p}\left(\int_{0}^{t_{1}}|A_{h,R}|^{\frac{m'}{q'}}\right)^{\frac{p}{\mu}}.
$$
Now we estimate the right hand side of \rife{sarca}: we first note that
$$
 \dys      \iakrt | f|\xi^p \leq \|f\|_{L^m (0,t_1; L^q (B_{R+\rho}))}  \left(\int_0^{t_1} |A_{k,R+\rho}|^{\frac{m'}{q'}}\right)^{\frac{1}{m' }}\,,
$$
and moreover
$$
\dys \int_0^{t_1} \int_{ A_{k,R}}  |T_j (\gku)|^{p_0} \leq c \|T_j (G_k(u))\|^{p_0}_{L^{d p_0} (Q^{t_1}_{R+\rho})} \left(\int_0^{t_1} |A_{k,R+\rho}|^{\frac{m'}{q'}}\right)^{\frac{1}{m' }}\,,
$$
where $d=\max \{ q,m\}$. Now, we observe that $f(t,x) \in L^{m_0} (0,t_1; L^{q_0} (B_{R+\rho}))$, $\forall m_0, q_0 $ such that
$1<m_0\leq m$, $1<q_0\leq q$. In particular we can choose $m_0, q_0$ such that
$$
\begin{array}{c}
m_0q=mq_0 \quad \mbox{and}\quad  \frac{1}{m_0} + \frac{N}{pq_0} = 1+\eps \,,\\[2.0 ex]
\dys \forall \eps < \min\left\{  \frac{m_0 q_0(N+p)+N(p-2)(q_0 (m_0-1)+m_0) }{p_0 d   }\,, \frac{N}{ p m_0}      \right\}\,.
\end{array}
$$
Using the first part of the Theorem we deduce that $u\in L^{\hat{s}} (Q_{R+\rho}^{t_1})$, where $\hat{s}= \frac{m_0q_0(N+p)+N(p-2)(q_0(m_0-1)+m_0) }{\eps}$. Since $\hat{s}\geq p_0 d$ we deduce
$$
\frac{c}{\rho }\dys \int_0^{t_1} \int_{ A_{k,R}}  |T_j (\gku)|^{p_0} \leq \frac{c_1 \|u\|_{L^{\hat{s}}(Q_{R+\rho}^{t_1})}}{\rho } \left(\int_0^{t_1} |A_{k,R+\rho}|^{\frac{m'}{q'}}\right)^{\frac{1}{m' }}\,.
$$
Gathering together the above informations, we finally deduce, using also that $\frac{\mu}{\sigma}=\frac{m'}{pq'}$ and passing to the limit with respect to $j$, that  there exists $C>0$ such that
$$
\int_{0}^{t_{1}}|A_{h,R} |^{\frac{m'}{q'}}\leq \frac{c}{(h-k)^{  \mu}\rho^\frac{\mu}{p}} \left(\int_{0}^{t_{1}}|A_{k,R+\rho} |^{\frac{m'}{q'}}\right)^{\frac{1}{m'}}\,.
$$
Since   \rife{326} is in force,  we have
$$
\frac{\mu}{m'p}=\frac{1}{q'}+\frac{p}{m' N}=1+\frac{p}{N}-\frac{p}{N}\left( \frac{1}{m} + \frac{N}{pq} \right)>1 \,,
$$
and so  we can apply Lemma \ref{stampa} to the function
$$
\zeta(h,d)=\int_{0}^{t_{1}}| A_{k,d}|^{\frac{m'}{q'}}  (t) \ dt\,.
$$
Thus the proof is complete for $0\leq t_{1}<T$.
As already remarked, it is clear that the choice of $t_{1}$ only depends on the data of the problem and thus we can iterate and conclude the same estimate in the whole cylinder in a finite number of steps.
\qed

\medskip

\noindent {\bf Acknowledgments. } 
The authors are grateful to CMA-University of Oslo and to CMUC-Departamento de Matem\'atica da Universidade de Coimbra for their kind hospitality.

\end{document}